\documentclass{article}

\newcommand{\Br}{\mathbb{R}}

\newcommand{\Prob}{\textnormal{Prob}}

\newcommand{\argmin}{\mathop{\rm argmin}}

\newcommand{\SCal}{\mathcal{S}}

\newcommand{\br}{\mathbb{R}}
\newcommand{\lmpz}{\lambda_{+}^0}
\newcommand{\lmmz}{\lambda_{-}^0}
\newcommand{\lmp}{\lambda_{+}}
\newcommand{\lmm}{\lambda_{-}}
\newcommand{\hsigma}{\hat{\sigma}}

\newcommand{\ba}{\begin{array}}
\newcommand{\ea}{\end{array}}

\usepackage{tcolorbox}
\usepackage{amsfonts,amscd,amssymb,amsmath,amsthm,enumerate}
\numberwithin{equation}{section}
\usepackage{array,mathtools,mathrsfs,bm}
\usepackage{threeparttable}
\usepackage{booktabs,caption,longtable,multicol,multirow}
\usepackage{subcaption,lscape}
\usepackage{makecell}
\usepackage{graphicx}
\usepackage{enumitem}
\setlist[enumerate]{noitemsep, topsep=2pt}
\setlist[itemize]{noitemsep, topsep=2pt}
\usepackage[colorlinks,citecolor=purple,linkcolor=blue]{hyperref}
\usepackage[nameinlink]{cleveref}
\usepackage[margin=1.5in]{geometry}
\usepackage[numbers]{natbib}
\usepackage{authblk}
\usepackage{algorithm,algorithmic}
\usepackage{pifont}
\usepackage{appendix}
\usepackage{multirow}
\usepackage{tabularx}
\usepackage{makecell}
\usepackage{eqparbox}

\usepackage{tikz}
\definecolor{ao(english)}{rgb}{0.0, 0.5, 0.0}
\definecolor{cadmiumgreen}{rgb}{0.0, 0.42, 0.24}
\definecolor{darkpastelgreen}{rgb}{0.01, 0.75, 0.24}

\usepackage{adjustbox}
\tikzset{%
  materia/.style={draw, fill=blue!20, text width=6.0em, text centered, minimum height=1.5em,drop shadow},
  etape/.style={materia, text width=8em, minimum width=10em, minimum height=3em, rounded corners, drop shadow},
  texto/.style={above, text width=6em, text centered},
  linepart/.style={draw, thick, color=black!50, -LaTeX, dashed},
  line/.style={draw, thick, color=black!50, -LaTeX},
  ur/.style={draw, text centered, minimum height=0.01em},
  back group/.style={fill=white!20,rounded corners, draw=black!50, dashed, inner xsep=15pt, inner ysep=10pt},
}

\tikzstyle{matheq} = [node distance=8.75cm, text width=21em, minimum width=1cm,
minimum height=2em, text centered,color=black!50]

\newtheorem{definition}{Definition}[section]
\newtheorem{proposition}{Proposition}[section]
\newtheorem{corollary}{Corollary}[section]
\newtheorem{theorem}{Theorem}[section]
\newtheorem{lemma}{Lemma}[section]
\newtheorem{remark}{Remark}[section]
\newtheorem{assumption}{Assumption}[section]

\title{Inexact and Implementable Accelerated Newton Proximal Extragradient Method for Convex Optimization}
\author[1]{\small Ziyu Huang}
\author[2]{\small Bo Jiang}
\author[2]{\small Yuntian Jiang\thanks{Correspondence to: yuntianjiang07@163.sufe.edu.cn}}

\affil[1]{\footnotesize School of Mathematical Sciences\\ Fudan University}
\affil[2]{\footnotesize School of Information Management and Engineering\\ Shanghai University of Finance and Economics}

\begin{document}
\maketitle

\begin{abstract}
    In this paper, we investigate the convergence behavior of the Accelerated Newton Proximal Extragradient (A-NPE) method\cite{monteiro2013accelerated} when employing inexact Hessian information. The exact A-NPE method was the pioneer near-optimal second-order approach, exhibiting an oracle complexity of $\Tilde{O}(\epsilon^{-2/7})$ for convex optimization. Despite its theoretical optimality, 
    there has been insufficient attention given to the study of its inexact version and efficient implementation. We introduce the inexact A-NPE method (IA-NPE), which is shown to maintain the near-optimal oracle complexity. In particular, we design a dynamic approach to 
    balance the computational cost of constructing the  Hessian matrix and the progress of the 
    convergence. Moreover,
    we show the robustness of the line-search procedure, which is a subroutine in IA-NPE, in the face of the inexactness of the Hessian. These nice properties enable the implementation of highly effective machine learning techniques like sub-sampling and various heuristics in the method. Extensive numerical results illustrate that IA-NPE compares favorably with state-of-the-art second-order methods, including Newton's method with cubic regularization and Trust-Region methods.
\end{abstract}
\section{Introduction}\label{Sec:Intro}
In this paper, we consider a generic unconstrained optimization problem  as follows:
\begin{equation}\label{Prob:main}
\min_{x\in\mathbb{R}^d} \ f(x):=g(x)+h(x),
\end{equation}
where $f$ is bounded below by $f^* > -\infty$, $g$ is convex and twice continuously differentiable with Lipschitz continuous Hessian, $h$ is convex and Lipshitz continuous but possibly non-differentiable.

From a theoretical standpoint, second-order methods are preferable for their ability to address ill-conditioned problems and exhibit a better oracle complexity. 
In the realm of convex optimization, several variations of prominent second-order methods stand out. Noteworthy examples include the trust-region method (TR) \cite{jiang2023universal}, the cubic regularized Newton method (CR) \cite{nesterov2006cubic}, and the gradient regularized Newton method (GR) \cite{doikov2022super,mishchenko2021regularized}.

Intriguingly, a historical review of accelerated second-order methods unveils their evolution over time. Nesterov \cite{nesterov2008accelerating} pioneered the field by proposing the initial accelerated second-order method, which improved the original version \cite{nesterov2006cubic}. This rate of acceleration was improved by Monteiro and Svaiter \cite{monteiro2013accelerated}, who introduced the A-NPE method, achieving an oracle complexity of $\tilde{O}(\epsilon^{-2/7})$. However, it's worth noting that their algorithm necessitates a search procedure in each iteration to determine the suitable step size, adding extra logarithmic complexity. In recent years, Arjevani et al. \cite{arjevani2019oracle} established a lower bound of $O(\epsilon^{-2/3p+1})$ for $p$-th order algorithms, highlighting that the A-NPE method happens to be nearly optimal up to a logarithmic factor. This observation has reignited interest in the A-NPE method. Simultaneously, researchers worldwide have independently extended the A-NPE method to accommodate higher-order information of objective functions. Three distinct groups\cite{gasnikov2019near} demonstrated that the modified method achieves an oracle complexity of $\Tilde{O}(\epsilon^{-2/3p+1})$, albeit still relying on a search procedure for the appropriate step size. Remarkably, a composite structure is allowed in the higher order A-NPE method \cite{jiang2021optimal}. 
Recently, Kovalev and Gasnikov \cite{kovalev2022first} and Carmon et al. \cite{carmon2022optimal} successfully eliminated the logarithmic factor in the complexity bound, making the A-NPE method a truly optimal algorithm. This marks a significant theoretical breakthrough in the field. 

Despite the notable theoretical achievements, it's crucial to acknowledge that the optimal complexity of the A-NPE method relies on exact Hessian information of the objective function, an impractical requirement in large-scale scenarios. In the realm of large-scale problems, the per-iteration computational complexity of second-order methods can become prohibitively expensive, primarily due to operations involving the Hessian. Addressing this challenge requires an examination of the convergence behavior when the Hessian is inexact. 

In contrast with standard second-order methods, the A-NPE method introduces an additional layer of complexity by requiring the identification of a suitable step size in every iteration, a parameter not known in advance. This necessitates a search procedure in each iteration,
and a good implementation of such a procedure is pivotal to the algorithm's performance.
To the best of our knowledge, most of the research on the A-NPE method focuses on the theoretical side, the only implementation by Carmon et al. \cite{carmon2022optimal} depends on exact Hessian information to derive the oracle complexity of $O(\epsilon^{-2/7})$, making it less practical for large-scale setting. (They also provide a CG routine to solve the subproblem, which fails to preserve the optimal oracle complexity.) This motivates us to analyze how the inexactness of the Hessian affects the search procedure. The challenge lies in proposing a robust algorithm that maintains theoretical complexity while being compatible with certain numerical heuristics.

Utilizing inexact Hessians is a widely employed technique to enhance the practical performance of second-order methods, demonstrating advantage from both theoretical and practical perspectives. The research focus is on designing algorithms using stochastic approximations to the Hessian. Various stochastic second-order methods have emerged, including but not limited to stochastic quasi-Newton methods \cite{byrd2016stochastic, schraudolph2007stochastic, wang2017stochastic}, stochastic cubic regularized Newton methods \cite{ masiha2022stochastic, tripuraneni2018stochastic, wang2019stochastic}, randomized cubic regularization methods \cite{doikov2018randomized}, stochastic trust-region methods \cite{blanchet2019convergence}, random subspace Newton methods \cite{gower2019rsn}, Hessian sketching methods \cite{ berahas2020investigation,hanzely2023sketch,lacotte2020effective,lacotte2021adaptive,pilanci2016iterative, pilanci2017newton}, sub-sampling methods \cite{agarwal2017second, bollapragada2019exact, byrd2011use, erdogdu2015convergence, kylasa2019gpu, li2020subsampled, liu2017inexact, roosta2019sub,xu2016sub,yao2021inexact}.

Among the extensive literature,  articles on inexact accelerated second-order methods are the most related to our study.
Ghadimi et al. \cite{ghadimi2017second} proposed an accelerated cubic regularized Newton method with an inexact Hessian, while Ye et al. \cite{ye2020nesterov} applied Nesterov's acceleration to enhance the approximate Newton method, both demonstrating favorable numerical performance. Song et al. \cite{song2019inexact} explored an accelerated inexact proximal cubic regularized Newton method with a complexity of $O(\epsilon^{-1/3})$ in the expectation sense. Chen et al. \cite{chen2022accelerating} and Kamzolov et al. \cite{kamzolov2023accelerated} investigated the accelerated adaptive cubic regularized Newton method, showcasing the use of sub-sampling and quasi-Newton methods to approximate the Hessian within this framework, both yielding a complexity of $O(\epsilon^{-1/3})$ and promising numerical results. Antonakopoulos et al. \cite{antonakopoulos2022extra} proposed a noise-adaptive accelerated second-order method with a universal global rate that adapts to the oracle's variance. Agafonov et al. \cite{agafonov2023inexact} introduced an accelerated inexact tensor method, demonstrating a complexity of $O(\epsilon^{-1/p+1})$ with access to the $p$-th order derivative. In recent work, Agafonov et al. \cite{agafonov2023advancing} proposed a second-order method with stochastic gradient and Hessian, proving its tight convergence bound
with respect to the variance of the gradient and the Hessian.
However, none of 
these inexact methods achieve near-optimal oracle complexity 
$\tilde{O}(\epsilon^{-2/7})$ for second-order methods.

In this paper, we propose the inexact Accelerated Newton Proximal Extra-gradient (IA-NPE) method, where an adaptive procedure is designed to determine the inexactness of the Hessian. In particular, we allow a relatively large error of the Hessian at the beginning of the algorithm
to reduce the cost of constructing the Hessian matrix. Note that the 
inexactness of the Hessian could affect the solution quality to the subproblem that appears in the main loop of the algorithm as well as the line search subroutine. 
Through meticulous analysis, we show that
the near-optimal oracle complexity of $\Tilde{O}(\epsilon^{-2/7})$ still holds for IA-NPE.
We discuss the adaptability of techniques such as Newton sketch and sub-sampling within the IA-NPE method. 
Furthermore, we illustrate that numerous heuristics can be seamlessly integrated into the method, enhancing the practical performance of the algorithm. We present extensive numerical experiments focused on logistic regression problems. The results showcase that the IA-NPE method performs comparably to mainstream algorithms in the context of machine learning problems.



The rest of the paper is organized as follows. In \autoref{Sec:Pre}, we introduce some basic definitions and assumptions required in the paper. In \autoref{Sec:framework}, we present the IA-NPE method in \autoref{alg.main alg} with its auxiliary search procedure in \autoref{alg.bisection}. In \autoref{Sec.ls}, we discuss how the inexactness of the Hessian information affects the auxiliary procedure and give a complexity analysis. In \autoref{Sec.convergence}, we analyze the overall complexity of the IA-NPE method. In \autoref{Sec.subroutine} and \autoref{Sec.numerical}, we demonstrate the applicability of modern machine-learning techniques to our algorithm and present promising results from numerical experiments.

\section{Preliminaries} \label{Sec:Pre}
In this section, we introduce the basic definitions and assumptions used in the paper. Denote the standard Euclidean norm in space $\mathbb{R}^d$ by \(\|\cdot\|\).  For a operator \(A : \mathbb{R}^d \rightarrow \Br^d\), its norm is defined as
\begin{equation*}
    \|A\| = \sup \{\|Ax\|:\|x\|\leq 1, \ x\in\br^d\}.
\end{equation*}
Throughout this paper, we refer to the following definition of $\epsilon$-optimality.
\begin{definition}
	 Given $\epsilon>0$, $x\in\br^d$ is said to be an $\epsilon$-optimal solution to problem~\eqref{Prob:main}, if
	\begin{equation}\label{result:optimality}
	f(x) - f^* \leq O(\epsilon),
	\end{equation}
\end{definition}
We can terminate the algorithm when we get points with such $\epsilon$-optimality.

Now we make some assumptions about the objective function.


\begin{assumption}
    \label{Assumption-Objective-Gradient-Hessian}
    The components $g$ and $h$ in \eqref{Prob:main} satisfy the following:
    \begin{itemize}
        \item $g$ and $h$ are proper closed convex functions.
        \item $g$ is twice continuously differentiable, the Hessian of $g$ is $L_2$-Lipschitz, i.e.,
        \begin{equation}
            \label{Def:Lipschitz-Hessian}
            \|\nabla ^2 g(x)-\nabla ^2 g(y)\| \leq L_2\|x-y\|, \ \forall x,y\in\br^d.
        \end{equation}
        \item $h$ is $L'$-Lipschitz continuous, i.e.,
        \begin{equation}
            \label{Def:Lipschitz-h}
            \vert h(x)-h(y)\vert \leq L'\|x-y\|, \ \forall x,y\in \br^d.
        \end{equation}
    \end{itemize}
\end{assumption}

Since our objective function concludes a non-differentiable part, we introduce the $\epsilon$-subdifferential of a proper closed convex function, whose basic properties are analyzed in Monteiro and Svaiter \cite[Section 2]{monteiro2013accelerated}.
\begin{definition}
For  $\epsilon \geq 0$, the $\epsilon$-subdifferential of a proper closed convex function $h:\br^d \to \Br \cup \{+\infty \}$ is the operator $\partial_{\epsilon}h:\br^d \rightrightarrows \Br^d$
\begin{equation}
\label{epsilon-subdifferential}
\partial_{\epsilon}h(x)=\left \{v\in \Br^d : h(y)\geq h(x)+\langle y-x,v\rangle-\epsilon, \ \forall y \in \br^d \right\} , \  x \in \br^d.
\end{equation}
\end{definition}

Regarding the smooth component, when the dimension of the problem is large, approximate Hessian is often used to reduce the computational cost. We consider the approximation as introduced in the following definition.

\begin{definition}\label{def.approximation}
The $\delta$-inexact second-order approximation of $g$ at $x$ is
\begin{equation}
\label{eq.approximationdelta}
g_{x,\delta}(y):=g({x})+\langle \nabla g({x}),y-{x}  \rangle +\frac{1}{2} \langle y-{x}, H\left(x\right)\left(y-x\right ) \rangle 
\end{equation}
where $H(\cdot)$ satisfies $\Vert H({x})-\nabla^2 g ({x}) \Vert <\delta$ and $H(x) \succeq 0$.
\end{definition}

Indeed, the above approximate Hessian can be constructed by various techniques, we will discuss this in \autoref{Sec.subroutine}.

\section{Overview of the IA-NPE Method}\label{Sec:framework}
\subsection{The IA-NPE Method and the Approximate Solution to the Subproblem}\label{framwork}
    \begin{algorithm}[!ht]
  \caption{The IA-NPE Method}\label{alg.main alg}
  \begin{algorithmic}
    \STATE{Initialization:} $x_0,y_0\in\br^d$, $A_0=0$, $C, \bar{\rho},\bar{\epsilon},\delta_{-1}>0,\delta_{\max}>\delta_{-1}$, $\gamma>1$, and $0<\sigma_l<\sigma_u<1,0<\hat{\sigma}<1$ with
$$\hat{\sigma}+\sigma_u <1,\quad \sigma_l(1+\hat{\sigma})<\sigma_u(1-\hat{\sigma}),\quad \tilde{\sigma}:=C+\sigma_u+\hat{\sigma} <1.$$ \\
    \FOR{$k = 0,1, 2, \cdots$}
    \STATE Set $\delta_k= \min\left\{\gamma\delta_{k-1},\delta_{\max}\right\}$, $\lambda_{k+1} = \frac{C}{\delta_k}$, construct the approximate Hessian $H(\cdot)$;\\
    \STATE Compute $(\Tilde{y}_{k+1},u_{k+1},\epsilon_{k+1})\in \operatorname{ANS}_{\hat{\sigma},\delta_k}(\lambda_{k+1},\Tilde{x}_k)$;\\
    \WHILE{$\lambda_{k+1}\|\tilde{y}_{k+1}-\tilde{x}_k\|< \frac{2\sigma_l}{L_2}$}
    \STATE Set $\delta_k = \delta_k/\gamma, \lambda_{k+1} = \frac{C}{\delta_k}$, update $H(\cdot)$ and $(\Tilde{y}_{k+1},u_{k+1},\epsilon_{k+1})\in \operatorname{ANS}_{\hat{\sigma},\delta_k}(\lambda_{k+1},\Tilde{x}_k)$;\\
    \ENDWHILE
    \IF{$\lambda_{k+1}\|\Tilde{y}_{k+1}-\tilde{x}_k\| >\frac{2\sigma_u}{L_2}$}
    \STATE Go to \textbf{bisection search} stage, compute $\lambda_{k+1}>0$ and $(\Tilde{y}_{k+1},u_{k+1},\epsilon_{k+1})\in \operatorname{ANS}_{\hat{\sigma},\delta_k}(\lambda_{k+1},\Tilde{x}_k)$ with
    \begin{equation}
        \label{eq.goalofsearch}
        \frac{2\sigma_l}{L_2}\le \lambda_{k+1}\|\tilde{y}_{k+1}-\tilde{x}_k\|\le \frac{2\sigma_u}{L_2},
    \end{equation}
    \ENDIF
    
    where
    \begin{equation}
    \label{eq.tildex}
    \begin{aligned}
        \tilde{x}_k &= \frac{A_k}{A_k+a_{k+1}}y_k+\frac{a_{k+1}}{A_k+a_{k+1}}x_k,\\
        a_{k+1}&= \frac{\lambda_{k+1}+\sqrt{\lambda_{k+1}^2+4\lambda_{k+1}A_k}}{2}
    \end{aligned}
    \end{equation}
    \STATE Set $y_{k+1}$ such that $f(y_{k+1}) \le f(\tilde{y}_{k+1})$ and 
\[v_{k+1}=\nabla g(\tilde{y}_{k+1})+u_{k+1}-\nabla g_{\tilde{x}_k,\delta}(\tilde{y}_{k+1}),\]
\[A_{k+1}=A_{k}+a_{k+1},\]
\[x_{k+1}=x_{k}-a_{k+1}v_{k+1}.\]
    \ENDFOR
  \end{algorithmic}
\end{algorithm}

In this section, we first present the IA-NPE method as in \autoref{alg.main alg} with its line-search subroutine as in \autoref{alg.bisection} to find the stepsize. Recall that the A-NPE method\cite{monteiro2013accelerated} exploits the second-order information of the smooth part $g$ to solve a proximal operator inexactly in each iteration. Different from the original version, we introduce an inner loop to determine the inexactness level $\delta$ and construct the approximate Hessian $H(\cdot)$ accordingly. Given the $\delta$-inexact second-order approximation of $g$, the current solution $x$ and stepsize $\lambda>0$, the IA-NPE method solve the following subproblem in each iteration
\begin{equation}\label{eq.second-order approximate ppa}
    y = \underset{u\in \Br^d}{\argmin} \ g_{x,\delta}(u)+h(u)+\frac{1}{2\lambda}\|u-x\|^2,
\end{equation}
and the optimality condition of \eqref{eq.second-order approximate ppa} is given by
\begin{equation}
    \label{eq.optimality condition ppa}
    v\in\left (\nabla g_{x,\delta}+\partial_\epsilon h\right )(y), \quad \lambda v+y-x=0.
\end{equation}
In \autoref{alg.main alg}, we solve \eqref{eq.second-order approximate ppa} inexactly and get the following approximate Newton solution.
\begin{definition}\label{Def.ANS}
Given $(\lambda,x) \in \Br_{++} \times \br^d$, error tolerance $\delta >0$ and $\hat{\sigma} \geq 0$, the triple $(y,u,\epsilon) \in \br^d \times \br^d \times \Br_+$ is called a $(\hat{\sigma},\delta)$-approximate Newton solution at $(\lambda,x)$ if
\begin{equation}\label{eq.ansu}
u \in (\nabla g_{x,\delta}+\partial_{\epsilon}h)(y),
\end{equation}
\begin{equation}\label{eq.ansstep}
\Vert \lambda u+y-x \Vert^2+2\lambda \epsilon \le \hat{\sigma}^2 \Vert  y-x \Vert^2.
\end{equation}
\end{definition} 
For simplicity, we denote the approximate solution in \autoref{Def.ANS}  for short as 
\begin{equation}\label{eq.ans for short}
(y,u,v)\in \operatorname{ANS}_{\hat{\sigma},\delta}(\lambda,x).
\end{equation}
The case $\delta=0$ in \autoref{Def.ANS} corresponds to a stricter approximate solution, which is based on the exact second-order information of the objective function.


\begin{definition}
    \label{def.AS}
    Given $(\lambda,x) \in \Br_{++} \times \br^d$, error tolerance $\tilde{\sigma} \geq 0$, the triple $(y,v,\epsilon) \in \br^d \times \br^d \times \Br_+$ is called a $\Tilde{\sigma}$-approximate solution at $(\lambda,x)$ if
\begin{equation}\label{eq.asu}
v \in (\nabla g_{x}+\partial_{\epsilon}h)(y),
\end{equation}
\begin{equation}\label{eq.asstep}
\Vert \lambda v+y-x \Vert^2+2\lambda \epsilon \le \tilde{\sigma}^2 \Vert  y-x \Vert^2,
\end{equation}
where $g_x(y)$ is the exact second-order expansion at $x$.
\end{definition}
We will frequently use the properties of the above approximate solution in the convergence analysis. Moreover,
we will later show that the two approximate solutions defined in \autoref{Def.ANS} and \autoref{def.AS} respectively are closely related.

Now we describe the line-search procedure in \autoref{alg.bisection}, where the subscript is omitted for simplicity. Note that as the way $\Tilde{x}$ is calculated in \autoref{alg.main alg}, $\Tilde{x}$ can be viewed as a continuous function of $\lambda$.
Though the Hessian in \autoref{alg.main alg} is inexact,  the bisection method as in \cite{monteiro2013accelerated} can still be adopted. 
Therefore, the problem considered in the line search procedure is as follows:\\
\textbf{Line-search Problem: }Given tolerance $\Tilde{\sigma}\geq 0$, $\Bar{\epsilon},\Bar{\rho}>0$, bounds $\alpha_+>\alpha_->0$, and a continuous curve $\Tilde{x}:[0,\infty)\rightarrow \br^d$ satisfies certain smoothness condition. The problem is to find a stepsize $\lambda>0$ and a $\tilde{\sigma}$-approximate solution $(y_\lambda,v_\lambda,\epsilon_\lambda)$ at $\left(\lambda,\Tilde{x}\left(\lambda\right)\right)$ such that 
\begin{equation*}
    \alpha_- \leq \lambda \|y_\lambda-\Tilde{x}(\lambda)\|\leq \alpha_+.
\end{equation*}
\begin{algorithm}[!ht]
  \caption{The Bisection Search}\label{alg.bisection}
  \begin{algorithmic}
    \STATE \textbf{Input:} Curve $\Tilde{x}:[0,+\infty)\rightarrow \Br$, right bracketing point $\lmpz>0$, the inexactness parameter $\delta,C$, and bounds $\alpha_+\ge \alpha_-\ge 0$ satisfying
    \begin{equation}
        \label{eq.boundforls}
        \alpha_-(1+\hat{\sigma}) < \alpha_+(1-\hat{\sigma})
    \end{equation}\\
    \STATE \textbf{Output:} stepsize $\lambda>0$ and $(y_\lambda,u_\lambda,\epsilon_\lambda)\in \operatorname{ANS}_{\hsigma,\delta}(\lambda,\Tilde{x}(\lambda))$ with $\alpha_-\leq \lambda\|y_\lambda-\Tilde{x}(\lambda)\|\leq \alpha_+$.
    \STATE \textbf{Bracketing Stage:} 
    \STATE compute $\lmmz$ satisfies \eqref{eq.left bracketing point} and set $x_-^0  = \Tilde{x}(\lmmz)$, compute $(y_-^0,u_-^0,\epsilon_-^0)\in \operatorname{ANS}_{\hsigma,\delta}(\lmmz,x_-^0)$;\\
    \STATE \textbf{Bisection Stage:} set $\lmm = \lmmz$, $\lmp = \lmpz$;\\
    \FOR{$k = 1,2,\cdots$}
    \STATE set $\lambda = \frac{\lmm+\lmp}{2}$, compute $(y_\lambda,u_\lambda,\epsilon_\lambda)\in \operatorname{ANS}_{\hsigma,\delta}(\lambda,\Tilde{x}(\lambda))$ and set $ v_{\lambda} = \nabla g(y_{\lambda}) - \nabla g_{x_{\lambda},\delta}(y_{\lambda})+u_{\lambda}$;
    \IF{$\lambda\|y_\lambda-\tilde{x}(\lambda)\|\in [\alpha_-,\alpha_+]$}
    \STATE \textbf{stop} and output the approximate solutions;\\
    \ELSE
    \STATE set 
    $$
    \left \{\begin{aligned}
        \lmp = \lambda, \ \lambda\|y_\lambda-\tilde{x}(\lambda)\|>\alpha_+,\\
        \lmm = \lambda, \ \lambda\|y_\lambda-\tilde{x}(\lambda)\|<\alpha_-.
    \end{aligned} \right .
    $$
    \ENDIF
    \ENDFOR
  \end{algorithmic}
\end{algorithm}

The dominating cost of \autoref{alg.main alg} and \autoref{alg.bisection} will be computing the approximate Newton solution, thus the complexity is evaluated in terms of the number of oracles required to compute such an approximate solution during the whole minimization process.

\begin{remark}
    \label{remark.early termination}
    If there exists a $\tilde{\sigma}$-approximate solution $(y,v,\epsilon)$ at $(\lambda,x)$ with $\epsilon\leq \Bar{\epsilon}$ and $\|v\| \leq \Bar{\rho}$, where $\Bar{\epsilon},\Bar{\rho}$ are given tolerance, then we can directly terminate the algorithm and output $y$ as an approximate solution of \eqref{Prob:main}. Therefore we always suppose there exists $\Bar{\epsilon},\Bar{\rho}\geq 0$ such that during the search procedure, the following holds for any $v$ and $\epsilon$:
    \begin{equation}
        \label{eq.relation v rho}
        \|v\| > \Bar{\rho}, \ \text{or } \epsilon> \Bar{\epsilon}.
    \end{equation}
\end{remark}

\subsection{Alternative Representation of Approximate Newton Solution}
Our analysis of the complexity of the algorithms 
depends on some existing results on the maximal monotone operator, thus we reformulate the approximate Newton solution and restate
some important definitions in the context of maximal monotone operators.
We put all proofs of the auxiliary technical results into the appendix for the coherence of the paper.

Note that as in \autoref{Assumption-Objective-Gradient-Hessian}, the minimization problem we consider can be viewed as a special case of the monotone inclusion problem with the following structure:
\begin{equation}\label{prob:main ls}
0 \in \mathcal{T}(x):=(\mathcal{G}+\mathcal{H})(x).
\end{equation}
Therefore, we make the following assumptions.
\begin{assumption}\label{assm.smoothness of mmo}
$\mathcal{G}: \br^d \to \br^d$ and $\mathcal{H}: \br^d \rightrightarrows \br^d$ satisfy the following:
\begin{itemize}
    \item $\mathcal{H}$ is a maximal monotone operator, also, $\|y\| \leq L'$ for any $x$ and $y\in \mathcal{H}(x)$.
    \item $\mathcal{G}$ is monotone and differentiable.
    \item $\mathcal{G}'$ is $L_2$-Lipschitz continuous on $\br^d$.
\end{itemize}
\end{assumption}

Similar to \autoref{def.AS}, we introduce the following approximate solutions of the proximal point iteration.
\begin{definition}
    \label{def.approximate solutions mmo}
    Given $(\lambda,x)\in\br_{++}\times \br^d$, tolerance $\tilde{\sigma}\geq 0$, the triple $(y,v,\epsilon)$ is said to be a $\tilde{\sigma}$-approximate solution at $(\lambda,x)$ if 
    \begin{equation}
        \label{eq.approximate solutions mmo}
        v\in \mathcal{T}^{\epsilon}(y), \ \|\lambda v+y-x\|^2+2\lambda\epsilon\leq \tilde{\sigma}^2\|y-x\|^2.
    \end{equation}
\end{definition}

Here we also permit an error tolerance $\delta$ for $\mathcal{G}'(\cdot)$, namely, we use an operator $\mathcal{P}(\cdot)$ with $\Vert \mathcal{G}'(x)-\mathcal{P} (x)\Vert <\delta$ at all possible $x$ during the search procedure to approximate $\mathcal{G}'(\cdot)$, we define the following inexact approximation and approximate Newton solution, corresponding to \autoref{def.approximation} and \autoref{Def.ANS}.

\begin{definition}\label{def.first-order approximation mmo}
For $x \in \br^d$, define the $\delta$-inexact first-order approximate of $\mathcal{T}_{x,\delta}:\br^d \to \br^d$ of $\mathcal{T}$ at $x$ as
\begin{equation}
\mathcal{T}_{x,\delta}(y)=\mathcal{G}_{x,\delta}(y)+\mathcal{H}(y), \ y \in \br^d,
\end{equation}
where $\mathcal{G}_{x,\delta}:\br^d \to \br^d$ is the $\delta$-inexact first-order approximate of $\mathcal{G}$ at $x$ given by:
\begin{equation}
\mathcal{G}_{x,\delta}(y)=\mathcal{G}(x)+\mathcal{P}(x)(y-x),
\end{equation}
where $\Vert \mathcal{P}(x)-\mathcal{G}'(x) \Vert <\delta$ and $\mathcal{P}$ should make $\mathcal{T}_{x,\delta}$ a maximal monotone operator.
\end{definition}

\begin{definition}\label{def.approximate newton solution mmo}
Given $(\lambda,x) \in \br_{++} \times \br^d$, error tolerance $\delta >0$ and $\hat{\sigma} \geq 0$, the triple $(y,u,\epsilon) \in \br^d \times \br^d \times \br_+$ is called a $(\hat{\sigma},\delta)$-approximate Newton solution at $(\lambda,x)$ if
\begin{equation}
u \in (\mathcal{G}_{x,\delta}+\mathcal{H}^{\epsilon})(y),
\end{equation}
\begin{equation}
\Vert \lambda u+y-x \Vert^2+2\lambda \epsilon \le \hat{\sigma}^2 \Vert  y-x \Vert^2,
\end{equation}
here $\mathcal{G}_{x,\delta}$ is the $\delta$-inexact first-order approximate of $\mathcal{G}$ at $x$ defined as in \autoref{def.first-order approximation mmo}.
\end{definition}
With a slight abuse of the notation, we still adopt the notion as \eqref{eq.ans for short} in the monotone operator setting. The next proposition shows that a $\Tilde{\sigma}$-approximate solution can be constructed with a $(\hat{\sigma},\delta)$-approximate Newton solution.
\begin{proposition}\label{prop:inexact approximate equation}
Let $(\lambda,x) \in \br_{++} \times \br^d$ and a $(\hat{\sigma},\delta)$-approximate Newton solution $(y,u,\epsilon)$ at $(\lambda,x)$ be given, and define $v:=\mathcal{G}(y)+u-\mathcal{G}_{x,\delta}(y)$.Then,
\begin{equation}\label{eq.inexact approxiamte equation v}
v \in (\mathcal{G}+\mathcal{H}^{\epsilon})(y) \subseteq \mathcal{T}^{\epsilon}(y),
\end{equation}
\begin{equation}\label{eq.inexact approxiamte equation lambda}
\Vert \lambda v+y-x \Vert^2+2\lambda \epsilon \le \left(\hat{\sigma}+\lambda\delta+\frac{L_2\lambda}{2} \left\Vert y-x \right\Vert\right)^2\Vert y-x \Vert^2,
\end{equation}
and
\begin{equation}\label{eq.bound v}
\Vert v \Vert \le \frac{1}{\lambda}\left(1+\hat{\sigma}+\lambda\delta+\frac{L_2\lambda}{2} \left\Vert y-x \right\Vert\right)\Vert y-x \Vert,
\end{equation}
\begin{equation}\label{eq.bound epsilon}
\epsilon \le \frac{\hat{\sigma}^2}{2\lambda} \Vert y-x \Vert^2.
\end{equation}
\end{proposition}

\section{Complexity of the Line-search Procedure}\label{Sec.ls}
In this section, motivated from the idea in Monteiro and Svaiter \cite[Section 7]{monteiro2013accelerated}, we will show that the complexity of the line-search procedure is logarithmic in the problem parameters, although the exact Hessian information is unavailable. 
\subsection{Preliminary Results}
For an approximate solution $(y,v,\epsilon)$ at $(\lambda,x)$ as in \autoref{def.approximate solutions mmo}, directly estimating the quantity $\lambda\|y-x\|$ may be difficult since the correspondence between $\lambda$ and the quantity is not single-valued. Thus it is necessary to find an approximation of the quantity that has a clear dependence of $\lambda$.

In order not to confuse the notations, consider a general maximal monotone operator $\mathcal{B}$, define for each $\lambda>0$,
\begin{equation}\label{prob.main.monotone operator}
y_\mathcal{B}(\lambda ;x):=(I+\lambda \mathcal{B})^{-1}(x), \  \varphi_\mathcal{B}(\lambda,x):=\lambda \Vert y_\mathcal{B}(\lambda ;x)-x \Vert.
\end{equation}
Where $y_\mathcal{B}(\lambda,x)$ is the exact proximal point iteration from $x$ with stepsize $\lambda$.

Now we introduce some basic properties of $\varphi_\mathcal{B}$ that will be needed in the analysis. It can be shown that $\varphi_\mathcal{B}(\lambda;x)$ can serve as a good approximation to $\lambda\|y-x\|$.

\begin{proposition}{(Monteiro and Svaiter \cite[Proposition 7.1]{monteiro2013accelerated})}\label{yinyong1}
For every $x \in \br^d$, the following statements hold:\\
(a)$\lambda >0 \to \varphi_\mathcal{B}(\lambda;x)$ is a continuous function;\\
(b)for every $0< \tilde{\lambda} \le \lambda$,
\begin{equation}
\frac{\lambda}{\tilde{\lambda}} \varphi_\mathcal{B}(\tilde{\lambda};x) \le \varphi_\mathcal{B}(\lambda;x) \le \left(\frac{\lambda}{\tilde{\lambda}}\right)^2 \varphi_\mathcal{B}(\tilde{\lambda};x)
\end{equation}
\end{proposition}
The following result shows that the quantity $\lambda\|y-x\|$ can be well-approximated by $\varphi_\mathcal{B}(\lambda;x)$. 
\begin{proposition}{(Monteiro and Svaiter \cite[Proposition 7.3]{monteiro2013accelerated})}\label{yinyong2}
Let $x \in \br^d$, $\lambda >0$ and $\tilde{\sigma}>0$ be given. If $(y,v,\epsilon)$ is a $\tilde{\sigma}$-approximate solution at $(\lambda,x)$,then
\begin{equation}
(1-\tilde{\sigma})\lambda \Vert y-x \Vert \le \varphi_\mathcal{B}(\lambda;x) \le (1+\tilde{\sigma}) \lambda \Vert y-x \Vert.
\end{equation}
\end{proposition}

\subsection{Analysis of The Bracketing Points}
The main goal of this section is to exploit the propositions of $\varphi_\mathcal{B}(\lambda;x)$ to find the bracketing points for the bisection procedure. Since we allow the derivative $\mathcal{G}'$ of the smooth part to be inexact, we have to make the following assumption about the approximation.

\begin{assumption}
The approximate operator $\mathcal{P}$ in \autoref{def.first-order approximation mmo} has the Lipschitz coefficient $L_2$.
\end{assumption}

By the mechanism of the bisection procedure and the fact that $\lmpz\delta = C$, we have the following observation.
\begin{proposition}\label{prop:approximation error}
During the search procedure, we always have $\lambda\delta \le C$.
\end{proposition}

The next proposition illustrates the behavior of $\varphi_{\mathcal{B}}$ in terms of $x$, which helps find the bracketing points.
\begin{proposition}\label{yuan7.5}
Given $\delta>0$ and let $x,\tilde{x} \in \br^d$ and $\lambda>0$ be given. Denote $\mathcal{B}:=\mathcal{T}_{x,\delta}$ and $\tilde{\mathcal{B}}:=\mathcal{T}_{\tilde{x},\delta}$, where  $\mathcal{T}_{x,\delta}$ and $\mathcal{T}_{\tilde{x},\delta}$ are the $\delta$-inexact first-order approximation of $\mathcal{T}$ at $x$ and $\Tilde{x}$. Then 
\begin{equation}\label{varphi-varphi}
\left \vert \varphi_\mathcal{B}\left(\lambda,x\right)-\varphi_{\tilde{\mathcal{B}}}\left(\lambda,\tilde{x}\right)\right\vert \le (1+2C)\lambda \Vert \tilde{x}-x \Vert+L_2\lambda^2 \Vert \tilde{x}-x \Vert^2+2L_2\lambda \Vert \tilde{x}-x \Vert \eta
\end{equation}
where
\begin{equation}
\eta:=\min\left \{\varphi_\mathcal{B}\left(\lambda;x\right),\varphi_{\tilde{\mathcal{B}}}\left(\lambda;\tilde{x}\right)\right\}.
\end{equation}
As a consequence, 
\begin{equation}
\varphi_\mathcal{B}(\lambda,x) \le (1+2C)\lambda \Vert \tilde{x}-x \Vert+L_2\lambda^2 \Vert \tilde{x}-x \Vert^2+\left(2L_2\lambda \left\Vert \tilde{x}-x\right \Vert +1\right)\varphi_{\Tilde{\mathcal{B}}}(\lambda;\Tilde{x}).
\end{equation}
\end{proposition}

With the help of the above propositions of $\varphi_{\mathcal{B}}$, we are now able to distinguish the bracketing points in the bisection procedure.

\begin{proposition}\label{prop.right bracketing point}
Let tolerance $\bar{v}>0$ and $\bar{\epsilon}>0$ , parameter $\hat{\sigma} \le 0$ and $\alpha>0$ be given. Then for any
\begin{equation}\label{eq.settolerance}
\lambda \geq \Lambda(\alpha):=\max \left\{ \sqrt{\frac{\alpha}{\bar{v}}\left(1+\hat{\sigma}+C+\frac{L_2\alpha}{2}\right)},  \left(\frac{\hat{\sigma}^2 \alpha^2}{2\bar{\epsilon}}\right)^{\frac{1}{3}} \right\},
\end{equation}
the $(\hat{\sigma},\delta)$-approximate Newton solution $(y,u,\epsilon)$ at $(\lambda, x)$ satisfies
$$\lambda \Vert y-x \Vert>\alpha$$
when $v = u -\mathcal{G}_{x,\delta}(y)+\mathcal{G}(y)$ such that \eqref{eq.relation v rho} holds.
\end{proposition}

\begin{proposition}\label{prop:min}
Let $(\lambda_+^0,x_+^0) \in \br_{++} \times \br^d$, and a $(\hat{\sigma},\delta)$-approximate Newton solution $(y_+^0,u_+^0,\epsilon_+^0)$ at $(\lambda_+^0,x_+^0)$ be given. Then, for any scalar $\alpha$ and $\lambda_-^0$ with
\begin{equation}\label{eq.lmpz}
0<\alpha \le \lambda_+^0 \Vert y_+^0-x_+^0 \Vert,
\end{equation}
\begin{align}\label{eq.lmmz}
&\lambda_-^0 \le \frac{\alpha (1-\hat{\sigma})\lambda_+^0}{(1+\hat{\sigma})(1+2L_2\theta_+^0)\lambda_+^0 \Vert y_+^0-x_+^0 \Vert+(1+2C)\theta_+^0+L_2(\theta_+^0)^2}, \\
\nonumber
&\theta_+^0:=\lambda_+^0 \Vert x_+^0-x_-^0 \Vert
\end{align}
the $(\hat{\sigma},\delta)$-approximate Newton solution $(y_-^0,u_-^0,\epsilon_-^0)$ at $(x_-^0,\lambda_-^0)$ satisfies
\begin{equation}\label{answer}
(1+\hat{\sigma})\lambda_-^0 \le (1-\hat{\sigma})\lambda_+^0,
\end{equation}
\begin{equation}\label{conclusion}
\lambda_-^0 \Vert y_-^0-x_-^0 \Vert \le \alpha.
\end{equation}
\end{proposition}

Note that in \eqref{eq.lmmz}, $\lambda_-^0$ depends on $x_-^0$, which is unknown at current stage. To get rid of such dependence, suppose that the curve $x(\cdot)$ additionally satisfies that 
\begin{equation}
\label{eq.curvecondition}
\Vert \tilde{x}(s)-\tilde{x}(t) \Vert \le \min\left\{\frac{M_0}{t} (s-t), M_1\left\Vert \tilde{x}(s)-\tilde{x}(0) \right\Vert\right\},  \   0 < t \leq s
\end{equation}
with some constant $M_0,M_1>0$, then we can set the left bracketing point as the following instruction.

\begin{corollary}\label{coro.left bracketing point}
In \autoref{alg.bisection}, if we set 
\begin{align}
    \label{eq.left bracketing point}
    &\lambda_{-}^0 =\frac{(1-\hat{\sigma})\alpha_-\lambda_{+}^0}{(1+\hat{\sigma})(1+2L_2\gamma_0)\lambda_{+}^0\Vert y_+^0-x_+^0 \Vert+(1+2C)\gamma_0+L_2(\gamma_0)^2}, \\
    \nonumber
    &\gamma_0 = M_1\lmpz \|x_+^0-x(0)\|.
\end{align}
then the following condition holds:
\begin{equation*}
\lambda_+^0 \Vert y_+^0-x_+^0\Vert \geq \alpha_+, \ \lambda_-^0 \Vert y_-^0-x_-^0\Vert \leq \alpha_-.
\end{equation*}
\end{corollary}

\subsection{Complexity of the Bisection Stage}
In this section, with the bracketing points at hand, we are able to analyze the complexity of the bisection stage of \autoref{alg.bisection}. We give the formal proof of \autoref{th:complexity_rough} here to show how the propositions contribute to the final analysis.



\begin{proposition}\label{prop:change}
Assume that $x_* \in \mathcal{T}^{-1}(0)=(\mathcal{G}+\mathcal{H})^{-1}(0)$, for any $\lambda>0$ and $x\in\br^d$, let $\mathcal{T}_{x,\delta}$ denotes the $\delta$-inexact first-order approximate of $\mathcal{T}$ at $x$, then there holds
\begin{equation}
\varphi_{\mathcal{T}_x}(\lambda;x)\le \frac{2L'\lambda^2}{1-C}+\frac{\lambda}{1-C} \Vert x-x_*\Vert+\frac{\lambda^2 L_2}{1-C} \Vert x-x_* \Vert^2,
\end{equation}
where $L'$ is the upper bound of $\mathcal{H}$ according to \eqref{prob:main ls}.
\end{proposition}

\begin{theorem}\label{th:complexity_rough}
The bisection stage of \autoref{alg.bisection} makes at most
\begin{equation}
2+\log\left [C_0^2M_0\lambda_+^0\left (\frac{1+2C+L_2 M_0 \lambda_+^0 +2\left(L_2+\frac{1}{M_0 \lambda_+^0}\right)(1+ \hat{\sigma})\alpha_-}{(1-\hat{\sigma})\alpha_+ - (1+\hat{\sigma})\alpha_-}\right)\right]
\end{equation}
oracle calls, where $\lambda_+^0$ is as in the definition in the framework above, and
\begin{equation}
\begin{aligned}
C_0=&\frac{(1+\hat{\sigma})(1+2L_2\gamma_0)\left[2L'\left(\lambda_+^0\right)^2+\lambda_+^0 d\left(x_+^0\right) + \left(\lambda_+^0\right)^2 L_2 d\left(x_+^0\right)^2\right]}{(1-C)(1-\hat{\sigma})^2\alpha_-}\\
&+\frac{\left[\left(1+2C\right)\gamma_0+L_2 \left(\gamma_0\right)^2\right](1-C)(1-\hat{\sigma})}{(1-C)(1-\hat{\sigma})^2\alpha_-},
\end{aligned}
\end{equation}
$\gamma_0$ is defined in \eqref{eq.left bracketing point} and $d(x_+^0)$ denotes the distance of $x_+^0$ to $x_*$.
\end{theorem}
{\it Proof}
First, we observe that the bracketing stage in \autoref{alg.bisection} makes one oracle call. By the mechanism of the bisection method, it follows that after $j$ bisection iterations,
\begin{equation}
\lambda_+-\lambda_-=\frac{1}{2^j}(\lambda_+^0-\lambda_-^0)\le \frac{1}{2^j}\lambda_+^0,
\end{equation}
and hence
\begin{equation}
j \le \log\left(\frac{\lambda_+^0}{\lambda_+-\lambda_-}\right).
\end{equation}
Assume now that the method doesn't stop at the $j$-th bisection iteration. Then, the values of $\lambda_-$ and $\lambda_+$ of this iteration satisfy
\begin{equation}
\lambda_+ \Vert y_{\lambda_+} -x(\lambda_+) \Vert >\alpha_+, \lambda_- \Vert y_{\lambda_-} -x(\lambda_-) \Vert <\alpha_-,
\end{equation}
and let $x_+:=\tilde{x}(\lambda_+)$, $x_-=\tilde{x}(\lambda_-)$. Let $\mathcal{B}_+:=\mathcal{T}_{x_+,\delta}$ and $\mathcal{B}_-:=\mathcal{T}_{x_-,\delta}$. Applying \autoref{yinyong2} twice, one time with $\mathcal{B}=\mathcal{B}_+, x=x_+$, and $(y,u,\epsilon)=(y_+,u_+,\epsilon_+)$, and the other with $\mathcal{B}=\mathcal{B}_-, x=x_-$, and $(y,u,\epsilon)=(y_-,u_-,\epsilon_-)$, we conclude that
\begin{equation}\label{thm5.4.1}
\varphi_+ :=\varphi_{\mathcal{B}_+}(\lambda_+;x_+)>(1-\hat{\sigma})\alpha_+, \ \varphi_- :=\varphi_{\mathcal{B}_-}(\lambda_-;x_-)<(1+\hat{\sigma})\alpha_-.
\end{equation}
On the other hand, it follows from \autoref{yinyong1} with $\mathcal{B}=\mathcal{B}_+, x=x_+, \tilde{\lambda}= \lambda_-$ and $\lambda=\lambda_+$, \autoref{yuan7.5} with $\lambda=\lambda_-, x=x_+$ and $\tilde{x}=x_-$, that
\begin{align}
\varphi_+ &=\varphi_{\mathcal{B}_+}(\lambda_+; x_+) \le \left(\frac{\lambda_+}{\lambda_-}\right)^2 \varphi_{\mathcal{B}_+}(\lambda_-; x_+)\\
& \le \left(\frac{\lambda_+}{\lambda_-}\right)^2\left[\left(1+2C\right)\theta+L_2\theta^2 +\left(1+2L_2\theta\right) \varphi_-\right],
\end{align}
from the smoothness condition \eqref{eq.curvecondition} we can conclude
\begin{equation}
\theta:=\lambda_- \Vert x_+-x_- \Vert \le M_0 (\lambda_+-\lambda_-),
\end{equation}
given the property assumed for the curve $x(\cdot)$. Hence, we conclude that
\begin{align*}
\varphi_+-\varphi_- \le & \left(\frac{\lambda_+}{\lambda_-}\right)^2 \theta (1+2C+L_2\theta +2L_2 \varphi_-)+ \left[\left(\frac{\lambda_+}{\lambda_-}\right)^2-1\right] \varphi_-\\
\le &\left(\frac{\lambda_+^0}{\lambda_-^0}\right)^2 M_0(1+2C+L_2 M_0 \lambda_+^0 +2L_2 \varphi_-)(\lambda_+-\lambda_-)\\
&+\frac{\lambda_+ + \lambda_-}{(\lambda_-)^2}\varphi_-(\lambda_+ - \lambda_-)\\
 \le& \left(\frac{\lambda_+^0}{\lambda_-^0}\right)^2 M_0(1+2C+L_2 M_0 \lambda_+^0 +2L_2 \varphi_-)(\lambda_+-\lambda_-)\\
 &+2\frac{\lambda_+^0}{(\lambda_-^0)^2}\varphi_-(\lambda_+ - \lambda_-)\\
 =&(\lambda_+ -\lambda_-)\left(\frac{\lambda_+^0}{\lambda_-^0}\right)^2 M_0 \left[1+2C+L_2 M_0 \lambda_+^0 +2\left(L_2+\frac{1}{M_0 \lambda_+^0}\right)\varphi_-\right].
\end{align*}
Combining the latter inequality with \eqref{thm5.4.1}, we then conclude that
\begin{equation}
\frac{1}{\lambda_+-\lambda_-} \le \left(\frac{\lambda_+^0}{\lambda_-^0}\right)^2 M_0 \left(\frac{1+2C+L_2 M_0 \lambda_+^0 +2\left(L_2+\frac{1}{M_0 \lambda_+^0}\right)\left(1+ \hat{\sigma}\right)\alpha_-}{\left(1-\hat{\sigma}\right)\alpha_+ - \left(1+\hat{\sigma}\right)\alpha_-}\right).
\end{equation}
We will now estimate the radio $\frac{\lambda_+^0}{\lambda_-^0}$. First note that \autoref{yinyong2} with $\lambda=\lambda_+^0, x=x_+^0, \mathcal{B}=\mathcal{T}_{x_+^0,\delta}$ and $(y,u,\epsilon)=(y_+^0, u_+^0, \epsilon_+^0)$, and \autoref{prop:change} with $\lambda=\lambda_+^0, x=x_+^0$ imply that
\begin{equation*}
\lambda_+^0 \Vert y_+^0 -x_+^0 \Vert \le \frac{\varphi_{\mathcal{T}_{x_+^0,\delta}}(\lambda_+^0;x_+^0)}{1-\hat{\sigma}} \le \frac{2L'(\lambda_+^0)^2+\lambda_+^0 d(x_+^0) + (\lambda_+^0)^2 L_2 d(x_+^0)^2}{(1-C)(1-\hat{\sigma})}.
\end{equation*}
Then using the latter inequality together with the definition of $\lambda_-^0$ in the framework, we can imply that
\begin{align*}
\frac{\lambda_+^0}{\lambda_-^0} =&\frac{(1+\hat{\sigma})(1+2L_2\gamma_0)\lambda_{+}^0\Vert y_+^0-x_+^0 \Vert+(1+2C)\gamma_0+L_2(\gamma_0)^2}{(1-\hat{\sigma})\alpha_-}\\
 \le &\frac{(1+\hat{\sigma})(1+2L_2\gamma_0)\left[2L'\left(\lambda_+^0\right)^2+\lambda_+^0 d\left(x_+^0\right) + \left(\lambda_+^0\right)^2 L_2 d\left(x_+^0\right)^2\right]}{(1-C)(1-\hat{\sigma})^2\alpha_-}\\
&+ \frac{+\left[\left(1+2C\right)\gamma_0+L_2\left(\gamma_0\right)^2\right](1-C)(1-\hat{\sigma})}{(1-C)(1-\hat{\sigma})^2\alpha_-}.
\end{align*}
Then we know
\begin{equation}
j \le \log\left[C_0^2M_0\lambda_+^0\left(\frac{1+2C+L_2 M_0 \lambda_+^0 +2\left(L_2+\frac{1}{M_0 \lambda_+^0}\right)(1+ \hat{\sigma})\frac{2\sigma_l}{L_2}}{(1-\hat{\sigma})\frac{2\sigma_u}{L_2} - (1+\hat{\sigma})\frac{2\sigma_l}{L_2}}\right)\right].
\end{equation}
The result now follows from the above inequality.
\qed
Furthermore, by ignoring the parameters, the above complexity can be simplified.
\begin{theorem}\label{thm.complexity_final}
\autoref{alg.bisection} makes at most
\begin{equation}\label{bound2}
O\left (\max\left \{\log L_2, \log d(x_+^0), \log M_1, \log M_0, \log L', \log \frac{1}{\Bar{\rho}}, \log \frac{1}{\Bar{\epsilon}} \right\}\right)
\end{equation}
oracle calls to compute a stepsize $\lambda_{k+1}>0$ that solves the line-search problem at $k$-th iteration.
\end{theorem}

\section{Complexity Analysis of the IA-NPE Method}\label{Sec.convergence}
\subsection{Convergence of the main loop in \autoref{alg.main alg}}
To establish the convergence of the IA-NPE method, the first thing is to check the complexity of the while loop of \autoref{alg.main alg}. Note that whenever $\lambda_{k+1}\|\Tilde{y}_{k+1}-\Tilde{x}_k\| \geq \frac{2\sigma_l}{L_2}$, the while loop will terminate, thus we come to the following lemma.
\begin{lemma}
    \label{lem.length while loop}
    In the $k$-th iteration of \autoref{alg.main alg}, the while loop will terminate after at most
    \begin{equation}
        \label{eq.length while loop}
        \log_{\gamma}\frac{\delta_{\max}\Lambda\left(\frac{2\sigma_l}{L_2}\right)}{C}
    \end{equation}
    oracle calls. Where $\Lambda(\cdot)$ is defined in \eqref{eq.settolerance}.
\end{lemma}
We now recall the following proposition that guarantees the iteration complexity of $O(\epsilon^{-2/7})$ for the main loop of the large step A-NPE method\cite{monteiro2013accelerated}.
\begin{proposition}
    \label{prop.ms proposition}
    There exists $0<\Tilde{\sigma}<1$ such that for every iteration $k$,
    \begin{equation}\label{con1}
    v_{k+1} \in (\nabla g+\partial_{\epsilon_{k+1}}h)(\tilde{y}_{k+1}) \subseteq \partial_{\epsilon_{k+1}}(g+h)(\tilde{y}_{k+1}),
    \end{equation}
    \begin{equation}\label{con2}
    \Vert \lambda_{k+1}v_{k+1}+\tilde{y}_{k+1}-\tilde{x}_k\Vert ^2+2\lambda_{k+1}\epsilon_{k+1}\le \Tilde{\sigma}^2 \Vert \tilde{y}_{k+1}-\tilde{x}_k\Vert^2,
    \end{equation}
\end{proposition}
It now remains to check if there exists a constant $\Tilde{\sigma}$ such that \autoref{prop.ms proposition} holds for \autoref{alg.main alg}. Luckily, when the error of the approximate Hessian is within $\delta_k$, 
such constant indeed exists and
the error brought by the inexact Hessians can be controlled.

\begin{lemma} \label{lem:rate}
In each iteration $k$, suppose that
\begin{equation}\label{eq.control the inexact}
\|\nabla ^2 g(\Tilde{x}_k)-H(\Tilde{x}_k)\|\leq \delta_k,
\end{equation}
we have
\begin{equation*}
v_{k+1} \in (\nabla g+\partial_{\epsilon_{k+1}}h)(\tilde{y}_{k+1})\subseteq \partial_{\epsilon_{k+1}}(g+h)(\tilde{y}_{k+1}),
\end{equation*}
\begin{equation}\label{lambdadelta}
\Vert \lambda_{k+1}v_{k+1}+\tilde{y}_{k+1}-\tilde{x}_k\Vert ^2+2\lambda_{k+1}\epsilon_{k+1}\le \Tilde{\sigma}^2 \Vert \tilde{y}_{k+1}-\tilde{x}_k\Vert^2,
\end{equation}
where $\Tilde{\sigma}=\sigma_u+\hat{\sigma} +C$.

\end{lemma}
The proof of \autoref{lem:rate} is similar to the proof of \autoref{prop:inexact approximate equation}, we omit it for simplicity. Now the convergence and the boundedness of the iterates can be guaranteed.
\begin{theorem}{(Monteiro and Svaiter \cite[Theorem 3.10, Theorem 6.4]{monteiro2013accelerated})}\label{thm.complexity framework}
Let $X^*$ be the set of optimal solutions and $x^*$ be the projection of $x_0$ onto $X^*$, denote $d_0 = \|x_0-x^*\|$, $\{x_k\},\{y_k\}$ are generated by \autoref{alg.main alg}, then we have 
\begin{gather*}
    \Vert x_k-x_*\Vert \le d_0, \ \Vert y_k-x^* \Vert \leq \left ( \frac{2}{\sqrt{1-\Tilde{\sigma}^2}}+1\right ) d_0,\\
    f(y_{k+1})-f_* \le \frac{3^{7/2}}{4\sqrt{2}} \frac{L_2 d_0^3}{\sigma_l \sqrt{1-\tilde{\sigma}^2}} \frac{1}{k^{7/2}}
\end{gather*}

\end{theorem}
\subsection{Total Complexity of the Algorithm}
Now we are able to analyze the total complexity of \autoref{alg.main alg}. We first note that the line-search stage in \autoref{alg.main alg} can be viewed as a special case of the problem we described in \autoref{Sec.ls}, it corresponds to the case where $\mathcal{H} = \partial h, \mathcal{G} = \nabla g$, the first-order approximation $\mathcal{T}_{x,\delta} = \nabla g_{x,\delta}+\partial h$. The curve $\Tilde{x}_k(\cdot)$ is calculated as in \eqref{eq.tildex}. We first show that the curve satisfies the smoothness condition \eqref{eq.curvecondition} and the subdifferentiable set is bounded.
\begin{lemma}{(Monteiro and Svaiter \cite[Lemma 7.13]{monteiro2013accelerated})}
    \label{lem.smoothness curve}
    The curve $\Tilde{x}_k(\cdot)$ satisfies \eqref{eq.curvecondition} with
    \begin{equation}
        \label{eq.curve parameter}
        M_0 = \left ( \frac{2}{\sqrt{1-\Tilde{\sigma}^2}}+1\right ) d_0, \ M_1 = 1.
    \end{equation}
\end{lemma}
\begin{lemma}\label{lemma:h}
If $h$ has a Lipschitz coefficient $L'$, then for any $\epsilon\geq0$, $x \in \br^d$ and any $v\in\partial h_\epsilon (x)$, we have
\begin{equation}\label{hL}
\Vert v \Vert \le L'.
\end{equation}
\end{lemma}

Now we can give the complete complexity result of \autoref{alg.main alg}.
\begin{theorem}
    \label{thm.complete complexity}
    If in each iteration $k$
    \begin{equation}
        \label{eq.assm total}
\|H(x_{k,\lambda})-\nabla ^2 g(x_{k,\lambda})\| \leq \delta_k,
    \end{equation}
    where $x_{k,\lambda}$ denotes all iterates during the line-search stage of the $k$-the iteration, then \autoref{alg.main alg} makes no more than 
    \begin{equation}
        \label{eq.complete complexity}
        O\left(L_2 d_0^3 \epsilon^{-2/7}\max \left\{\left(\log L_2\right )^2,\log d_0, \log L',\left(\log \left (1/\Bar{\rho}\right)\right )^2,\left (\log \left (1/\Bar{\epsilon}\right)\right )^2\right\}\right)
    \end{equation}
    oracle calls to find an approximate solution $y_k$ satisfies
    \begin{equation*}
        f(y_k)-f^*\leq \epsilon
    \end{equation*}
    or a $\Tilde{\sigma}$-approximate solution $(y_k,v_k,\epsilon_k)$ with 
    \begin{equation*}
        \|v_k\| \leq \Bar{\rho}, \ \epsilon_k\leq \Bar{\epsilon}.
    \end{equation*}
\end{theorem}
{\it Proof}
    This theorem is a direct consequence of \autoref{thm.complexity_final}, \autoref{lem.length while loop}, \autoref{lem.smoothness curve}, \autoref{thm.complexity framework}.
\qed

\section{Subroutines for Approximating the Hessian}\label{Sec.subroutine}
In this section, we show that when the objective function in \eqref{Prob:main} has certain special structures, some techniques can be applied to approximate the Hessian of the smooth part. For example, when the dimension of the problem is high, we can use a subspace approximation to the Hessian. In addition, we can apply a sub-sampling technique when the objective function has a finite-sum structure. For the sub-sampling technique, we provide a formal analysis and show that the technique is compatible with the IA-NPE algorithm both theoretically and practically.

\subsection{Subspace Approximation}
Note that in \autoref{alg.main alg}, we first determine the error of the inexact Hessian $\delta_k$ and then construct the approximate Hessian, when $\delta_k$ is relatively large, we can use the subspace approximation or the Newton sketch technique. 

Take the Newton sketch technique as an example, when the square root of the Hessian of the objective function is easily computable, we can apply this technique. Let's say an objective function with a structure of
\begin{equation*}
    f(x) = \sum_{i=1}^n f_i(a_i^T x),
\end{equation*}
the square root of the Hessian is $\nabla ^2 f(x)^{1/2}= \textbf{diag}\left(f_i''\left(a_i^T x\right)^{1/2}\right) A$, where $A$ denotes the data matrix.

Specifically, at every iteration, a random sketch matrix $S_k\in \br^{m\times d}$ is defined, where $m<d$. The approximate Hessian was defined to be 
\begin{equation}
    \label{eq.sketch hessian}
    H(x) := \nabla^2 g(x)^{1/2} S_k^T S_k \nabla^2 g(x).
\end{equation}
Theoretically, it is difficult to establish the relation between the error $\delta_k$ and the dimension $m$. However, we can dynamically increase the sketch dimension when the current error is seemingly larger than the pre-given inexactness $\delta_k$.

An inner loop to approximate the Hessian by sketching can be designed. We can start with a small dimension, construct the sketch Hessian, and check whether the left and right bracketing points are well-defined, specifically, we should check if
\begin{equation*}
    \lambda_+^0 \Vert y_+^0-x_+^0\Vert \geq \alpha_+, \ \lambda_-^0 \Vert y_-^0-x_-^0\Vert \leq \alpha_-
\end{equation*}
holds in the bisection stage, if not,  we can double the sketching dimension, i.e., we set $m\to 2m$ and construct the sketch Hessian again. Empirically, when the dimension approaches the effective dimension as in \cite{lacotte2021adaptive}, which is usually less than $d$, we may ensure that we have a good approximate. Even in the worst case, when the dimension approaches $d$, the sketching Hessian will be a good approximation. Thus the inner loop at most has a logarithmic complexity.

\subsection{Sub-Sampling Approximation} \label{Section2:Sampling}
When objective function has the finite-sum structure
\begin{equation*}
    f(x) = \sum_{i=1}^n f_i(x),
\end{equation*}
the sub-sampling technique can help reduce the computational cost, we assume each component $f_i$ satisfies the following assumption.

\begin{assumption}\label{SSAssumption-Objective-Gradient-Hessian}
Each $f_j(x)$ has a composite form $f_j(x)=g_j(x)+h_j(x)$, where the following conditions are assumed to hold:
\begin{itemize}
\item The objective function $f(x):=\frac{1}{n} \sum_{j=1}^n f_j(x)$ is convex, $g(x):=\frac{1}{n} \sum_{j=1}^n g_j(x)$ and $h(x):=\frac{1}{n} \sum_{j=1}^n h_j(x)$ are both proper closed convex functions.
\item  $g_j$ is twice continuously differentiable with its gradient and Hessian being Lipschitz continuous, i.e., there are $0<L_{1,g_j}, L_{2,g_j}<\infty$ such that for any $x, y\in\br^d$
\begin{align}
\left\| \nabla g_j(x) - \nabla g_j(y)\right\| \leq &  \ L_{1,g_j} \left\|x-y\right\|, \label{SSDef:Lipschitz-Gradient} \\
\left\| \nabla^2 g_j(x) - \nabla^2 g_j(y)\right\| \leq &  \ L_{2,g_j} \left\|x-y\right\|. \label{SSDef:Lipschitz-Hessian}
\end{align}
\item $h_j$ is Lipschitz continuous, i.e., there exists $0<L_{h_j}<\infty$ such that for any $x, y\in\br^d$
\begin{align}
\left\| h_j(x) - h_j(y)\right\| \leq &  \ L_{h_j} \left\|x-y\right\|. \label{SSDef:Lipschitz-h}
\end{align}
\end{itemize}
\end{assumption}

In the rest of the section, we define 
\begin{equation*}
L_1=\max_{1\le j\le n} L_{1,g_j} , \  L_2=\max_{1\le j \le n} L_{2,g_j} ,\ L'=\max_{1\le j\le n} L_{h_j} .
\end{equation*}

The following technical lemmas are mostly from Xu et al. \cite{xu2020newton} so we omit the proof here.

Let $\SCal$ and $\left|\SCal\right|$ denote the sample collection and its cardinality, $p_i=\Prob\left(\xi=i\right)$ denotes the probability that the index $i$ is chosen, and define 
\begin{equation}\label{Condition: SSHessian}
\Tilde{H}(x) = \frac{1}{n\left|\SCal\right|}\sum_{j\in\SCal} \frac{1}{p_j}\nabla^2 g_j(x).
\end{equation}
When $n$ is very large, such random sampling can significantly reduce the per-iteration computational cost as $\left|\SCal\right| \ll n$.

One natural sampling scheme is to use the uniform sampling, i.e., $p_i=\frac{1}{n}$ for $i=1,2,\ldots,n$. There are other sampling schemes, we refer interested readers to \cite{roosta2019sub,xu2016sub,xu2020newton}.

The following lemmas reveal how many samples are required to get a sub-sampled Hessian within a given accuracy if
the indices are sampled uniformly with replacement. 


\begin{lemma}\label{cor:uniform}
Suppose \autoref{Assumption-Objective-Gradient-Hessian} holds and for $x\in\br^d$,
$\Tilde{H}(x)$ is constructed from \eqref{Condition: SSHessian} with $p_j=\frac{1}{n}$. When sample size
\begin{equation}
\label{eq.samplesizeuni}
\left|\SCal\right| \geq \frac{16L_1^2}{\kappa^2} \log\left(\frac{2Nd}{\delta'}\right)
\end{equation} 
for given $0<\kappa, \delta' <1$. Then for any $N$ points $x_i \in \br^d, i=1,\ldots,N$ we have
\begin{equation*}
\Prob\left( \bigcup_{i=1}^N \left\{\left\| \Tilde{H}(x_i) - \nabla^2 g(x_i)\right\| \geq \kappa \right\} \right) < \delta'.
\end{equation*}
\end{lemma}

In view \autoref{Assumption-Objective-Gradient-Hessian} and the way we construct the sub-sampled Hessian matrices stated in \eqref{Condition: SSHessian}, we have the following remark.

\begin{remark}
The sub-sampled Hessian $\Tilde{H}(x_k)$ constructed by \eqref{Condition: SSHessian} with sample size satisfying \eqref{eq.samplesizeuni} has the Lipschitz coefficient $L_2$.
\end{remark}

The sub-sampled Hessian won't be changed during the search procedure to preserve the Lipschitz continuity. We will specify the way of constructing the approximate Hessian $H(\cdot)$ in \autoref{alg.main alg}.
\begin{lemma}
    \label{lem.construct approximate hessian ss}
    Given the overall failure probability $\delta_0$, in the $k$-th iteration, if we set 
\begin{equation}
    \label{eq.set sub-sampled Hessian}
    \delta' = O(\delta_0 \epsilon^{2/7}), \ \kappa = \frac{\delta_k}{2}, \ N = \max \left\{\log L_2,\log d_0, \log L',\log \left (1/\Bar{\rho}\right),\log \left (1/\Bar{\epsilon}\right)\right\}
\end{equation}
in \eqref{eq.samplesizeuni} and set the approximate Hessian as 
\begin{equation}
    \label{eq.good approximate hessian ss}
    H(x) = \Tilde{H}(x)+\frac{\delta_k}{2} I,
\end{equation}
then 
\begin{equation}
    \label{eq.approximate hessian ss}
    \| H(x_{k,\lambda})-\nabla^2 g(x_{k,\lambda})\|\leq \delta_k, \ H(x)\succeq 0
\end{equation}
with probability $1-O(\delta_0 \epsilon^{2/7})$.
\end{lemma}

To describe the convergence behavior of \autoref{alg.main alg} in terms of probability, we have the following equivalent result.
\begin{theorem}
    \label{thm.prob convergence}
    If we construct the approximate Hessian as instructed in \autoref{lem.construct approximate hessian ss}, then when \autoref{alg.main alg} makes 
    $$
    O\left(L_2 d_0^3 \epsilon^{-2/7}\max \left\{\left(\log L_2\right )^2,\log d_0, \log L',\left(\log \left (1/\Bar{\rho}\right)\right )^2,\left (\log \left (1/\Bar{\epsilon}\right)\right )^2\right\}\right)
    $$
    oracle calls, with probability $1-\delta_0$ we have an approximate solution $y_T$ with
    \begin{equation*}
        f(y_T)-f^*\le O(\epsilon)
    \end{equation*}
    or we have a $\Tilde{\sigma}$-approximate solution $(y_T,v_T,\epsilon_T)$ with 
    \begin{equation*}
        \|v_T\| \leq \Bar{\rho}, \ \epsilon_T\leq \Bar{\epsilon}.
    \end{equation*}
\end{theorem}

\section{Numerical Experiment}\label{Sec.numerical}
In this section, we present the results of our numerical experiment. During the experiment, we compare our algorithm with state-of-the-art algorithms in the context of regularized logistic regression problems. The findings indicate that our algorithm is well-suited for tackling large-scale machine-learning problems. These experiments were carried out on a Lenovo laptop equipped with a 3.10 GHz CPU and 16GB of memory.

\textbf{Problem.}
The model of regularized logistic regression is given by
\begin{equation}
\label{eq.logreg}
\underset{x\in\mathbb{R}^d}{\min}\quad f(x)=\frac{1}{n}\sum_{i=1}^n\log(1+e^{-b_i w_i^T x})+\frac{\alpha}{2}\| x\|,
\end{equation}
where $\{(w_i,b_i)\}_{i=1}^n$ are the problem data with $w_i$ and $b_i\in\{-1,1\}$ denote the characteristic and label separately, the regularizer $\alpha$ is set as 1e-5.

We perform experiment on the following six data sets from LIBVIM$\footnote{https://www.csie.ntu.edu.tw/~cjlin/libsvmtools/data sets/}$, their statistics are summarized in the following Table \ref{dataset}:
\begin{table}[h]
\caption{The Statistics of six data sets}
\centering
\label{dataset}
\begin{tabular}{|c|c|c|}
\hline
Name          & Instance No.($n$) & Feature No($d$). \\ \hline
phishing       & 11055      & 68          \\ \hline
a8a     & 22696        & 123         \\ \hline
a9a & 32561      & 123           \\ \hline
ijcnn1 & 49900  & 22                        \\ \hline
SUSY          &5,000,000     & 18          \\ \hline
HIGGS         &11,000,000	 & 28          \\ \hline
\end{tabular}
\end{table}

\textbf{Heuristic Adaptive Scheme.}
The line-search procedure is the key to the implementation, we propose a heuristic adaptive search method. In the $k$-th iteration, the sample size is chosen inversely proportional to the square norm of the gradient and is bounded below and above by some constants, which depend on the statistic of the data set. Motivated by Carmon et al. \cite{carmon2022optimal}, we set the initial stepsize based on the value we obtained in the previous loop, the theoretical analysis in \autoref{Sec.ls} ensures that the stepsize will always lie in the bracketing interval.

In \eqref{eq.logreg}, the nonsmooth part $h(x)=0$, it is the following two inequalities that guarantee the convergence of \autoref{alg.main alg}:
    \begin{align}
    \label{eq.convergencegua1}
    \lambda_{k+1}\| \tilde{y}_{k+1}-\tilde{x}_{k+1}\|&\geq \frac{2\sigma_l}{L_2}\\
    \label{eq.convergencegua2}
    \| \lambda_{k+1} \nabla f(\tilde{y}_{k+1})+\tilde{y}_{k+1}-\tilde{x}_{k}\| &\leq \Tilde{\sigma} \| \tilde{y}_{k+1}-\tilde{x}_{k} \|
    \end{align}

In our implementation, we do not guarantee the above two inequalities strictly. We set a threshold to evaluate whether the current line-search procedure is costly and adjust the parameters to satisfy it dynamically. In the $k$-th iteration, for \eqref{eq.convergencegua1}, when the current search appears to be simple, i.e., the length of the search procedure does not exceed the threshold, we will set the initial stepsize in the next iteration larger, say $\lambda_{k+1}=2\lambda_k$, to pursue more aggressive performance. For \eqref{eq.convergencegua2}, we make the parameter $\Tilde{\sigma}$ an adaptive parameter with lower and upper bounds $\underline{\sigma}<\bar{\sigma}\in\left(0,1\right)$, we denote it as $\sigma_k$, when the line-search procedure is simple, we decrease $\sigma_k$ to achieve more aggressive performance, if the line-search procedure seems to be too difficult, we increase $\sigma_k$ and redo the line-search at the current point.


As for the acceleration technique, despite its nice theoretical properties, it does not always 
show superiority in the local scenario. The A-NPE framework can be seen as a method to accelerate the gradient regularized Newton method, during the implementation,  we switch to the sub-sampled gradient regularized Newton method when the acceleration period is almost finished. To clarify, in our implementation, after 40 accelerations and when the current iterate seems no better than the previous iteration, i.e., $\frac{\vert f(y_{k+1})-f(y_k)\vert}{\vert f(y_k) \vert} \leq 0.1$, we set $\lambda_k = \| \nabla f(y_k)\| ^{\frac{3}{2}}$ and update
$$y_{k+1}= \underset{x \in \br^d}{\arg\min} \ \langle \nabla f\left(y_k\right),x-y_k\rangle +\frac{1}{2}\langle x-y_k, H\left(y_k\right)\left(x-y_k\right)\rangle+\frac{\lambda_k}{2}\|x-y_k\|^2.$$ 

\textbf{Experiment Setting.}
We set the adaptive parameters used in the Algorithm as $\sigma_0=0.9,\underline{\sigma}= 0.7,\bar{\sigma}=0.95$, and the threshold is set as 2. The final stopping criterion is set as $\|\nabla f(y_k)\| < 10^{-7}$. The initial point of is set as a random point of Gaussian distribution with zero mean and $10^5$ covariance. We denote \autoref{alg.main alg} with the above heuristics as adaptive A-NPE method(AANPE).

\subsection{Comparison with State-of-art Deterministic Algorithms}
In our experiment, to show the effectiveness of the sub-sampling technique, we compare our algorithm with the following deterministic algorithms, including the classic trust-region method(TR)\cite{conn2000trust}, the Broyden-Fletcher-Goldfarb-Shanno(BFGS) method and its limited memory version(LBFGS)\cite{nocedal1999numerical}, the adaptive cubic regularized Newton method(ARC)\cite{cartis2011adaptive1,cartis2011adaptive2}, and the adaptive accelerating cubic regularized Newton method (AARC)\cite{jiang2020unified}. 

We adopt the implementation of ARC and TR in the public package\cite{kohler2017sub} $\footnote{https://github.com/dalab/subsampled\_cubic\_regularization}$ with the
default parameters and the full batch taken. To solve the ARC and TR subproblem, we use the Lanczos method and GLTR method separately. For BFGS and LBFGS, we adopted the well-tuned algorithms from the public package 'scipy' with full-batch gradient information used.

\begin{figure}[!ht]

  \begin{minipage}[b]{0.49\textwidth}
    \includegraphics[width=\textwidth]{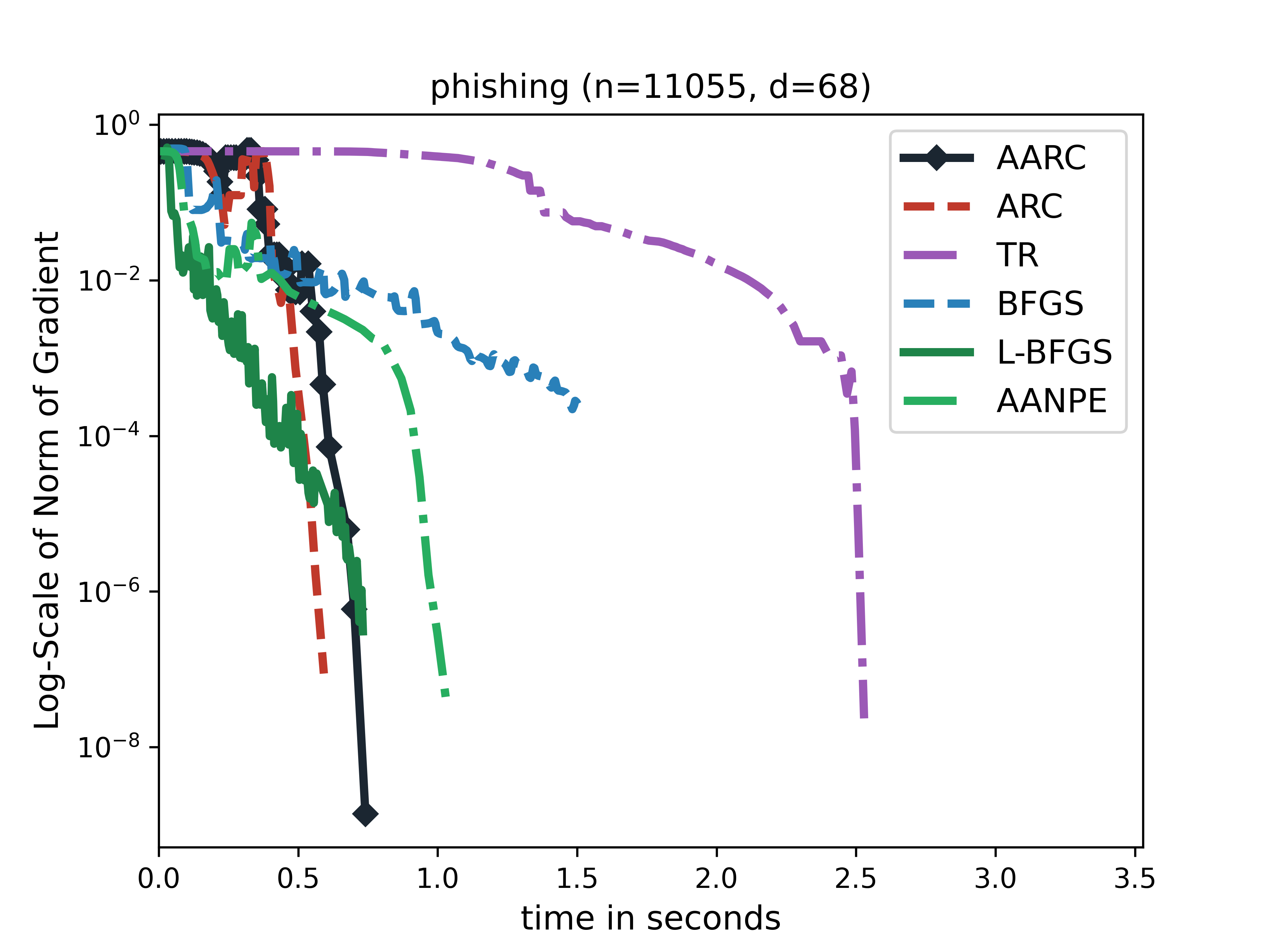}
  \end{minipage}
  \begin{minipage}[b]{0.49\textwidth}
    \includegraphics[width=\textwidth]{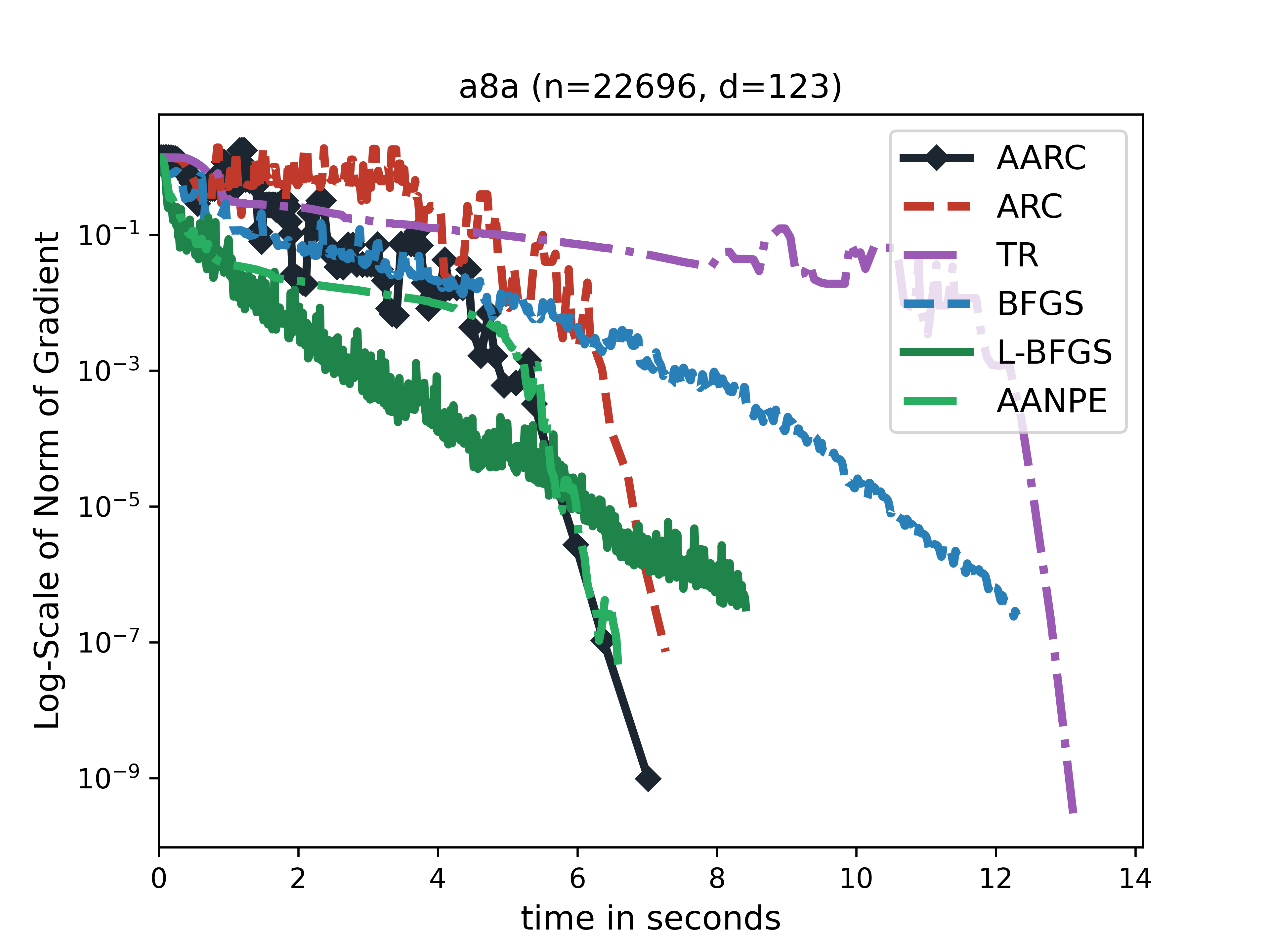}
  \end{minipage}
  \begin{minipage}[b]{0.49\textwidth}
    \includegraphics[width=\textwidth]{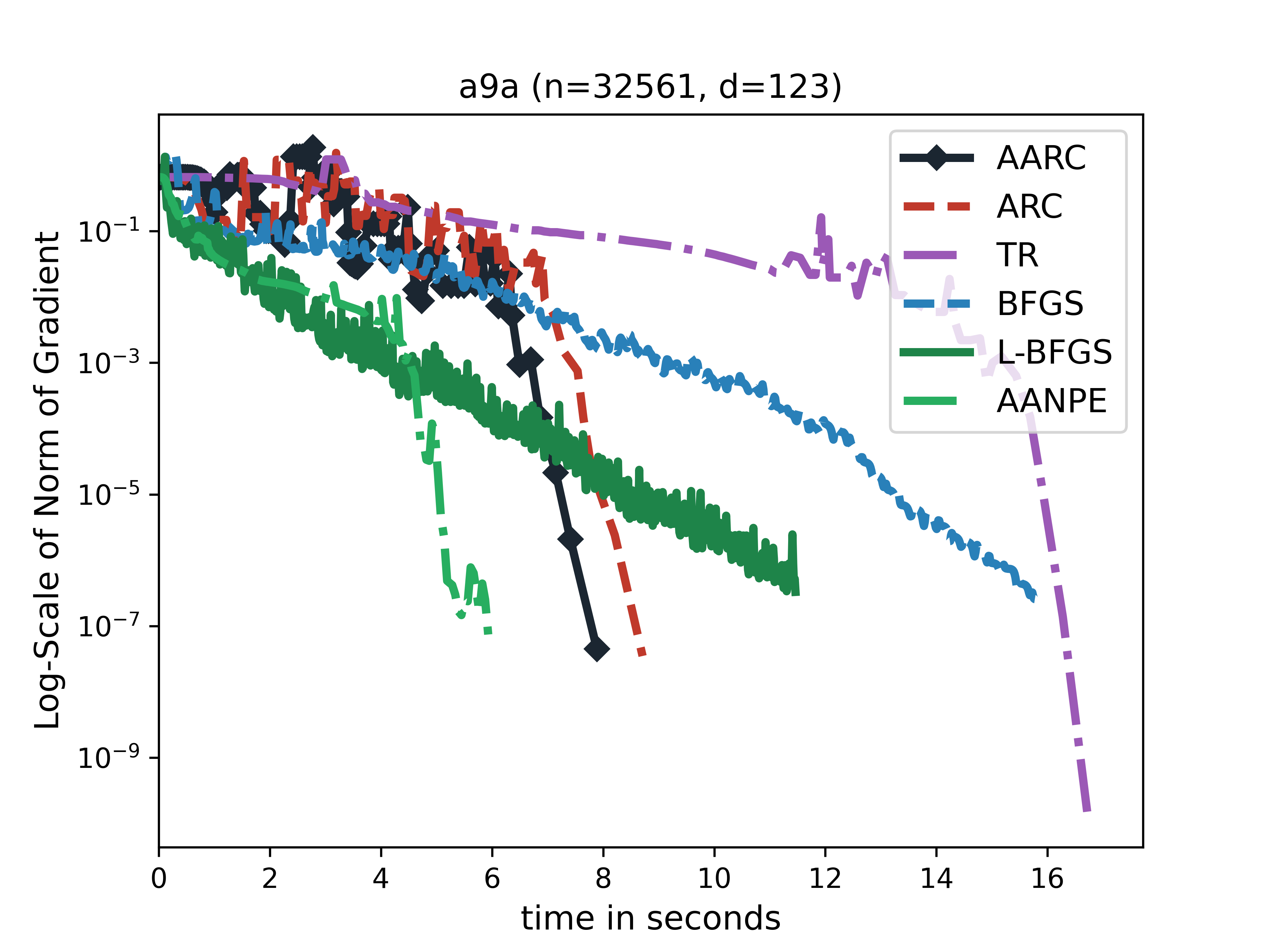}
  \end{minipage}
  \begin{minipage}[b]{0.49\textwidth}
    \includegraphics[width=\textwidth]{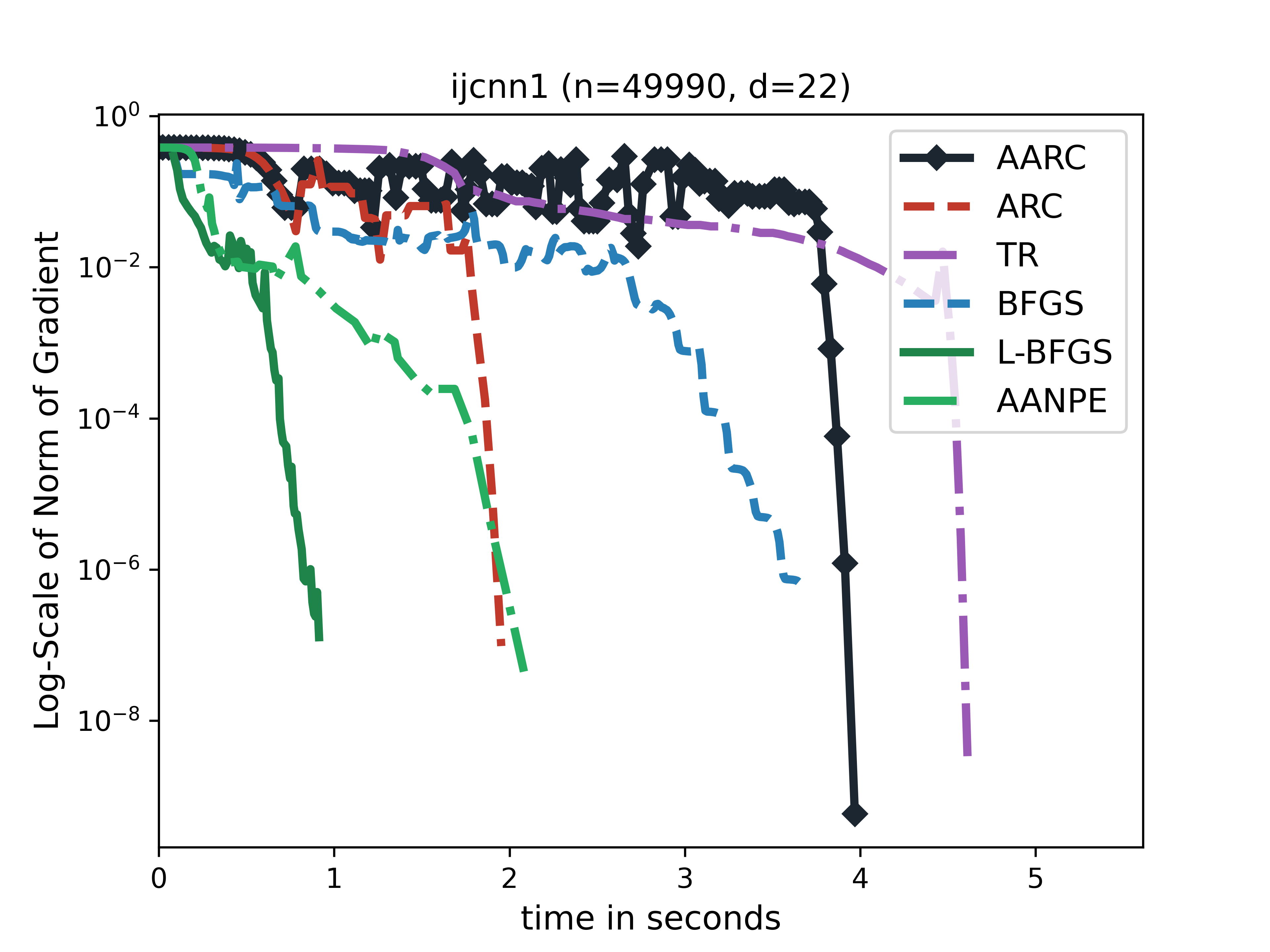}
  \end{minipage}
  \begin{minipage}[b]{0.49\textwidth}
    \includegraphics[width=\textwidth]{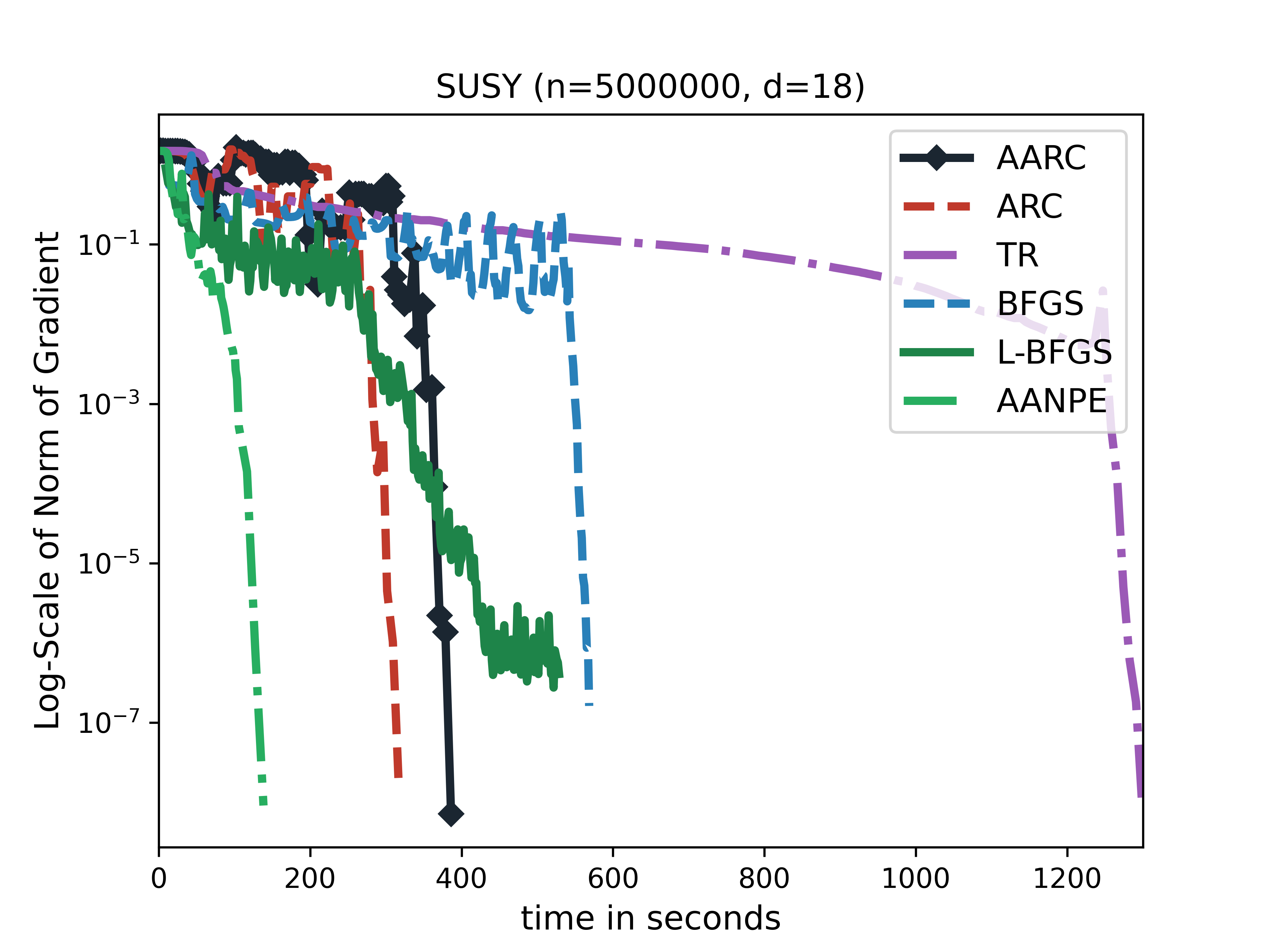}
  \end{minipage}
  \begin{minipage}[b]{0.49\textwidth}
    \includegraphics[width=\textwidth]{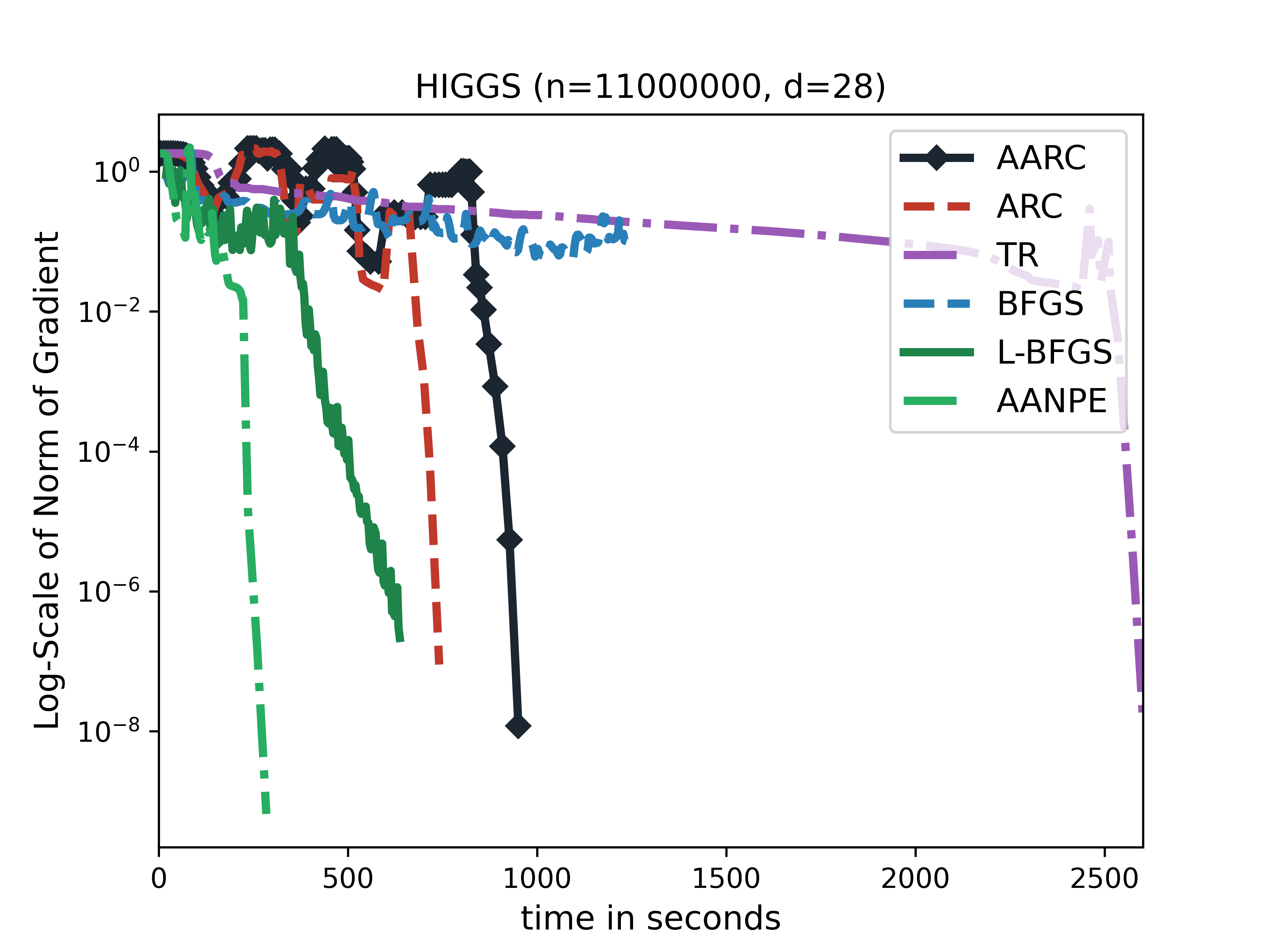}
  \end{minipage}
  \caption{
    Log-scale norm of gradient v.s. time.
  }
  \label{fig.gravstime}
\end{figure}

\begin{figure}[!ht]

  \begin{minipage}[b]{0.49\textwidth}
    \includegraphics[width=\textwidth]{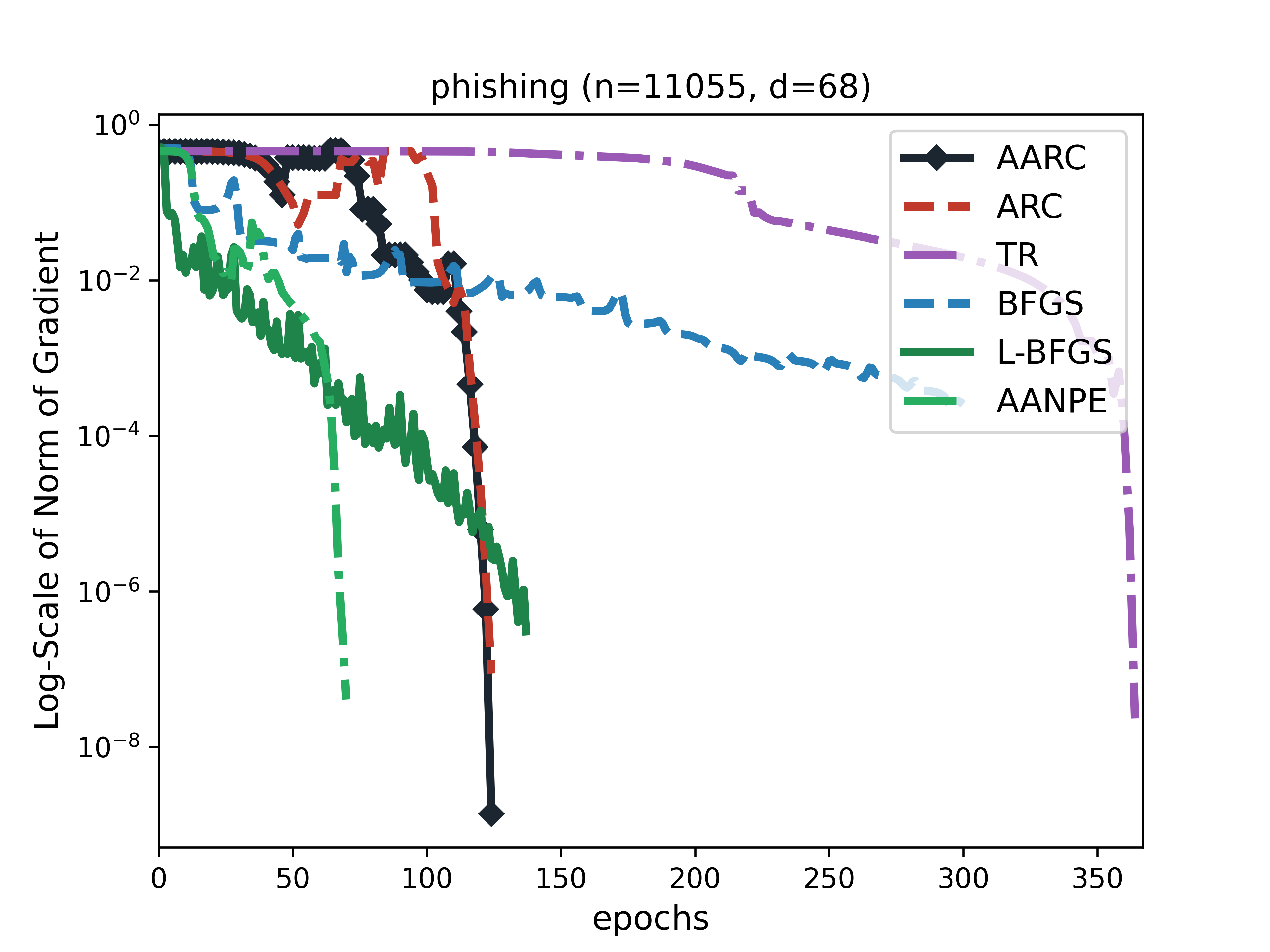}
  \end{minipage}
  \begin{minipage}[b]{0.49\textwidth}
    \includegraphics[width=\textwidth]{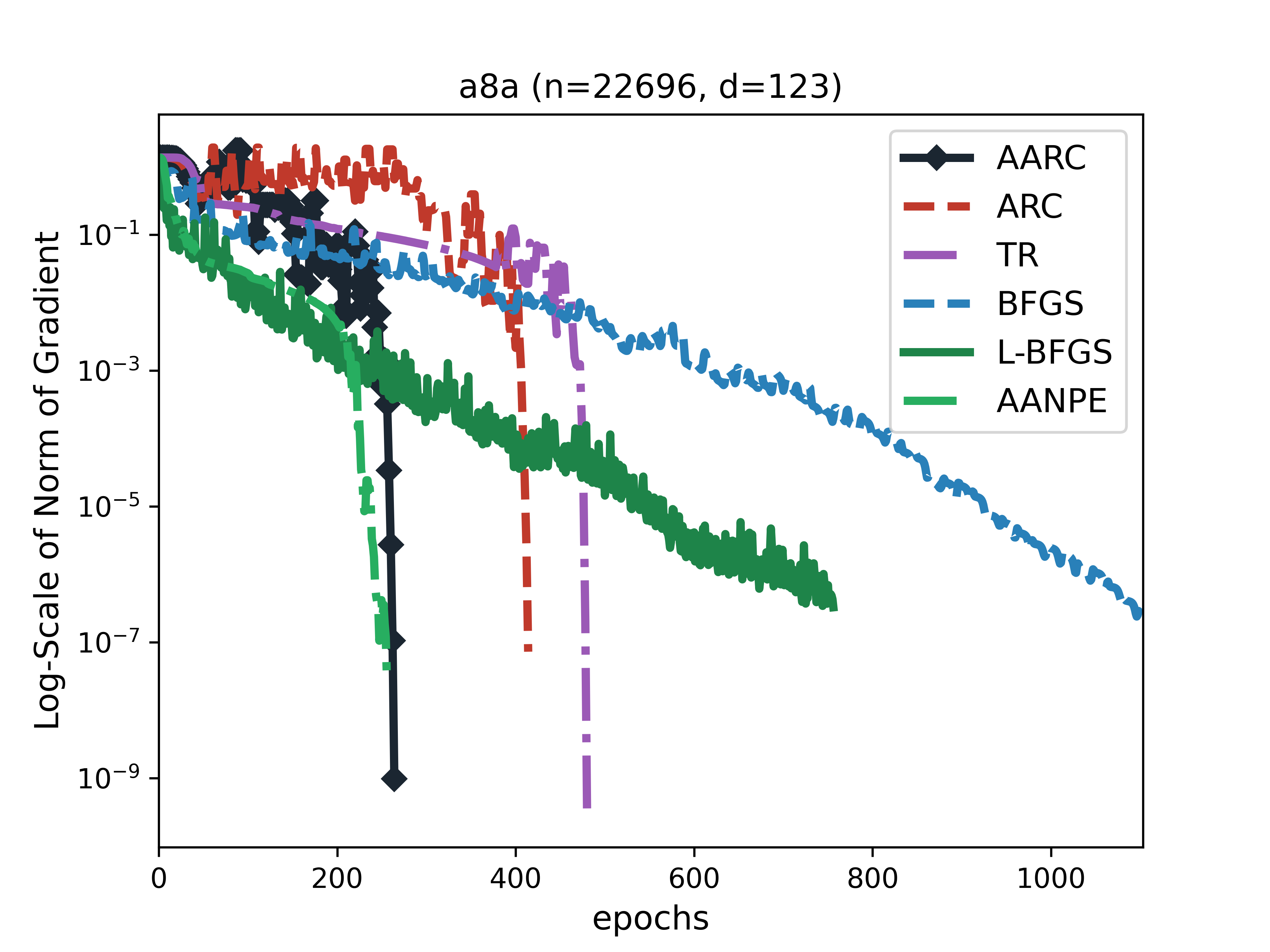}
  \end{minipage}
    \begin{minipage}[b]{0.49\textwidth}
    \includegraphics[width=\textwidth]{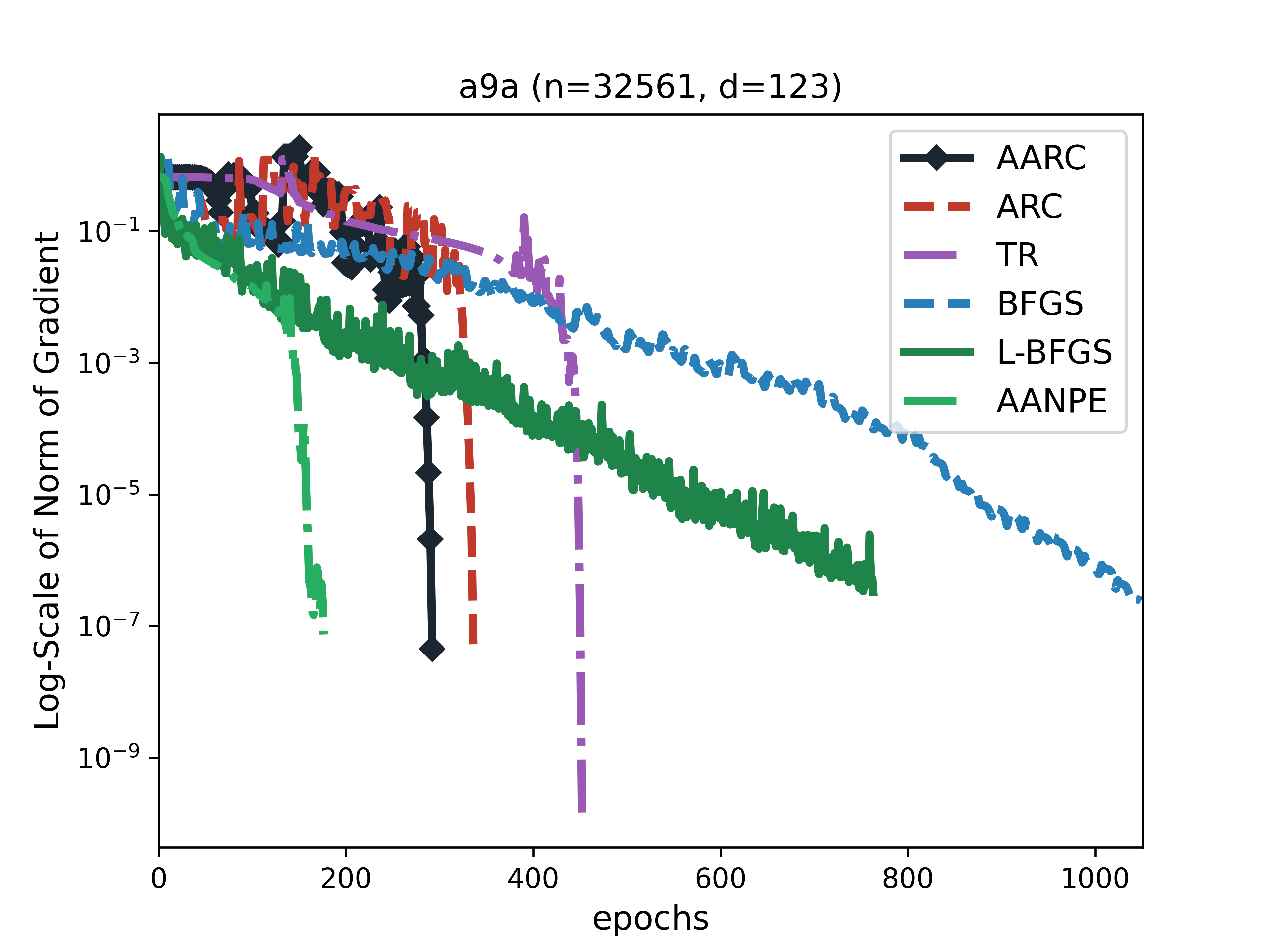}
  \end{minipage}
  \begin{minipage}[b]{0.49\textwidth}
    \includegraphics[width=\textwidth]{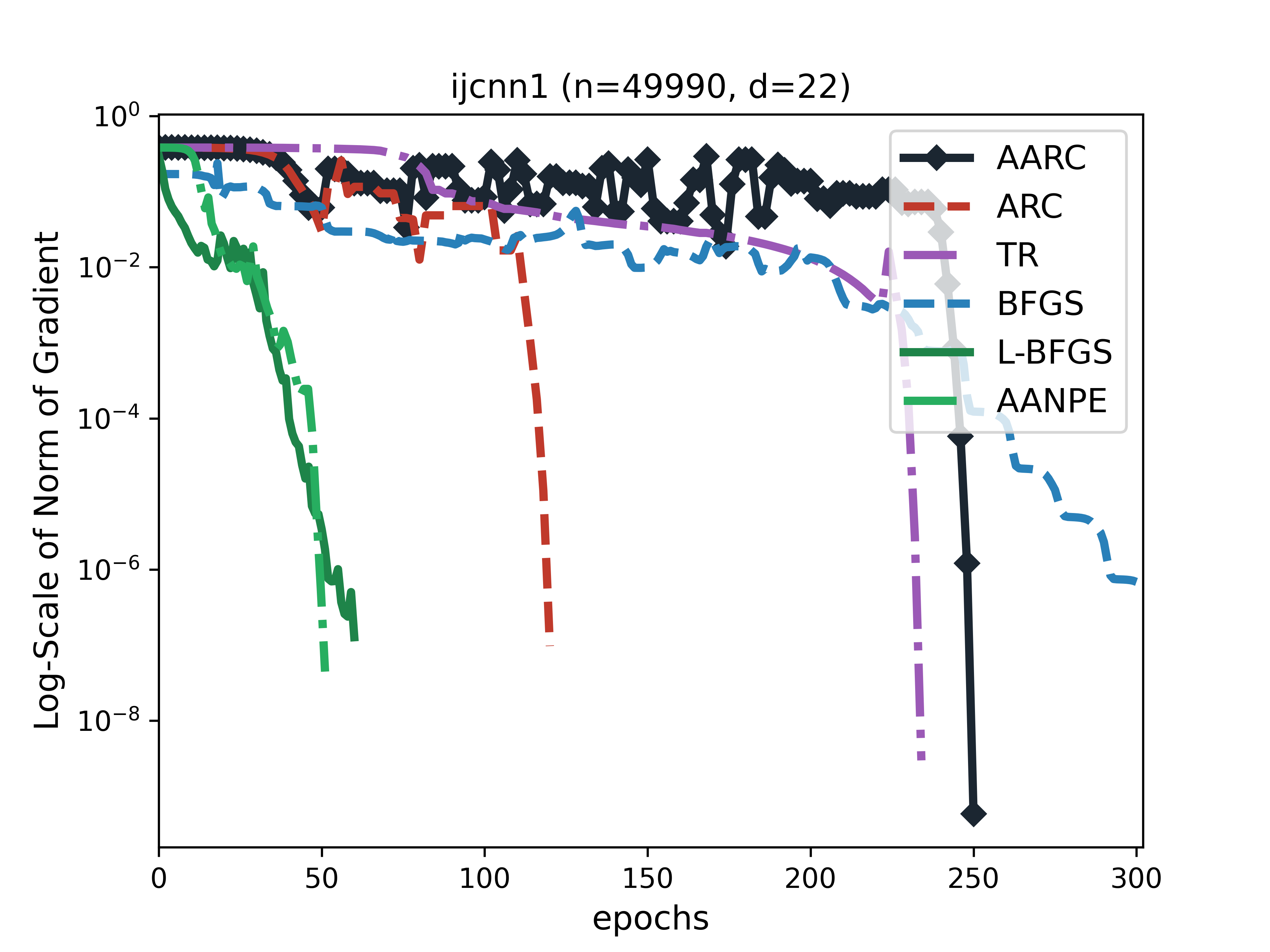}
  \end{minipage}
  \begin{minipage}[b]{0.49\textwidth}
    \includegraphics[width=\textwidth]{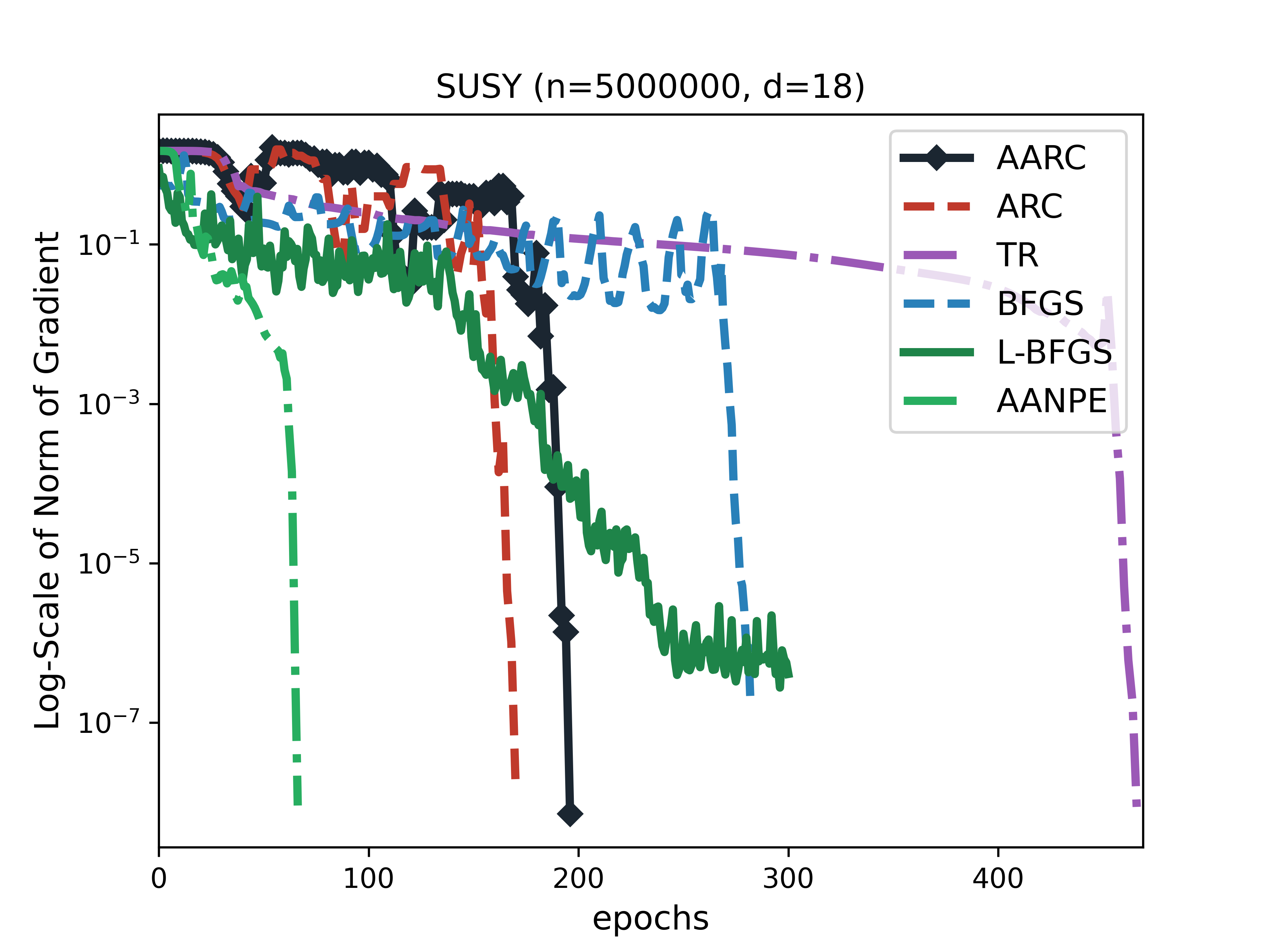}
  \end{minipage}
  \begin{minipage}[b]{0.49\textwidth}
    \includegraphics[width=\textwidth]{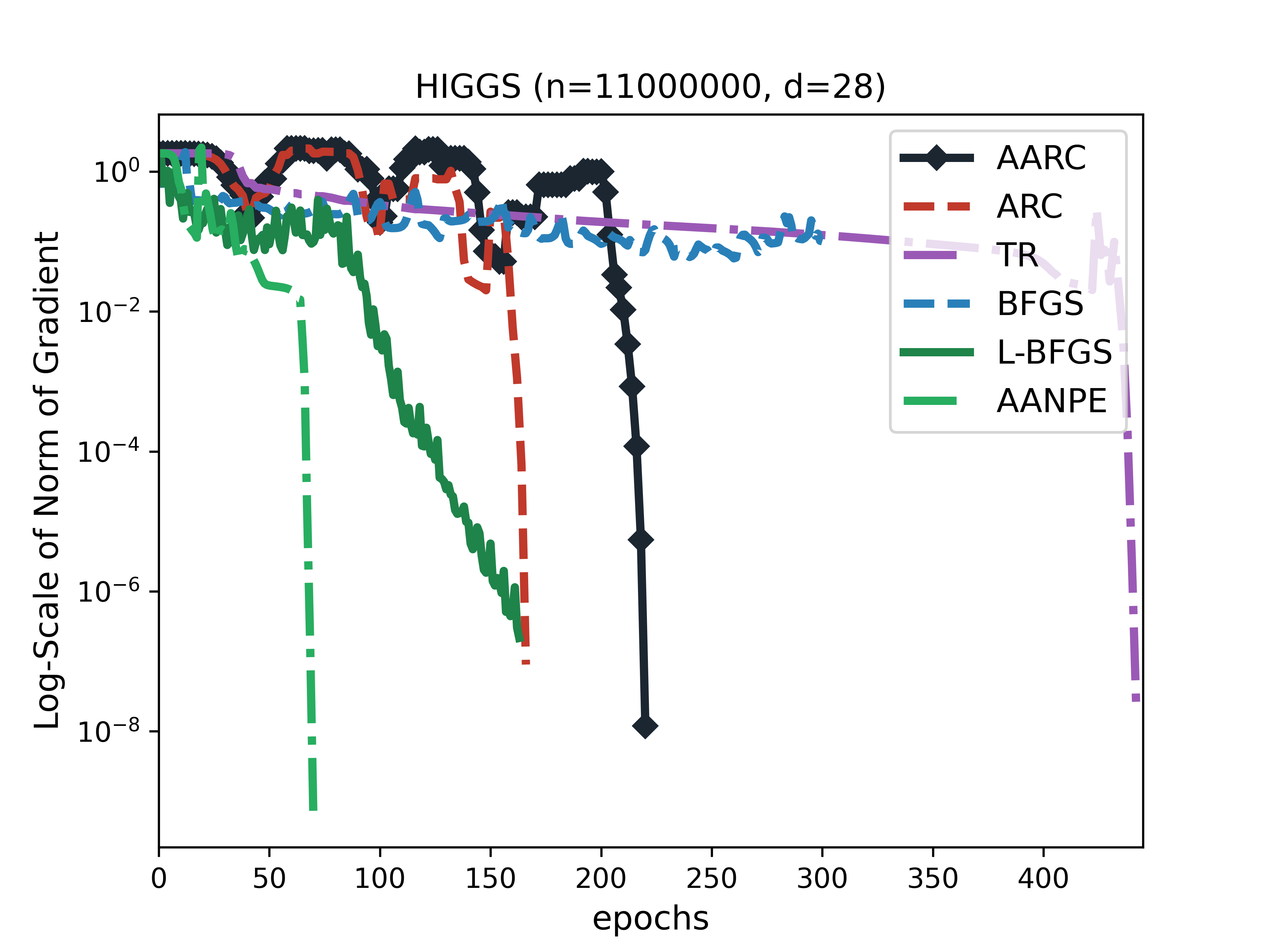}
  \end{minipage}
  \caption{
    Log-scale norm of gradient v.s. No. epochs.
  }
  \label{fig.gravsiter}
\end{figure}

The results are presented in \autoref{fig.gravstime} and \autoref{fig.gravsiter}. It shows that our algorithm is comparable to the state-of-art algorithms. Especially, on data sets with a large number of instances such as 'SUSY' and 'HIGGS', our algorithm outperforms all deterministic algorithms significantly, this phenomenon shows that the the sub-sampling technique indeed accelerates the algorithm.

\subsection{Comparison with Two Different Types of Sub-sampled Algorithms}
In the second experiment, we compare our algorithm with two other sub-sampled second-order algorithms, one is the sub-sampled cubic regularized Newton method(SCR)\cite{kohler2017sub}, the other one is its accelerated variant, the sub-sampled accelerating adaptive cubic regularized Newton method(SACR)\cite{chen2022accelerating}. In the implementation, we set the initial point as zero vector and the number of the accelerated iterations as 5. Regarding the sampling scheme for the two other algorithms, we use the adaptive sampling scheme equipped with the default parameters in SCR as in their code. For SACR, We also follow their instruction in \cite{chen2022accelerating} on choosing the parameters. 

\begin{figure}[!ht]

  \begin{minipage}[b]{0.49\textwidth}
    \includegraphics[width=\textwidth]{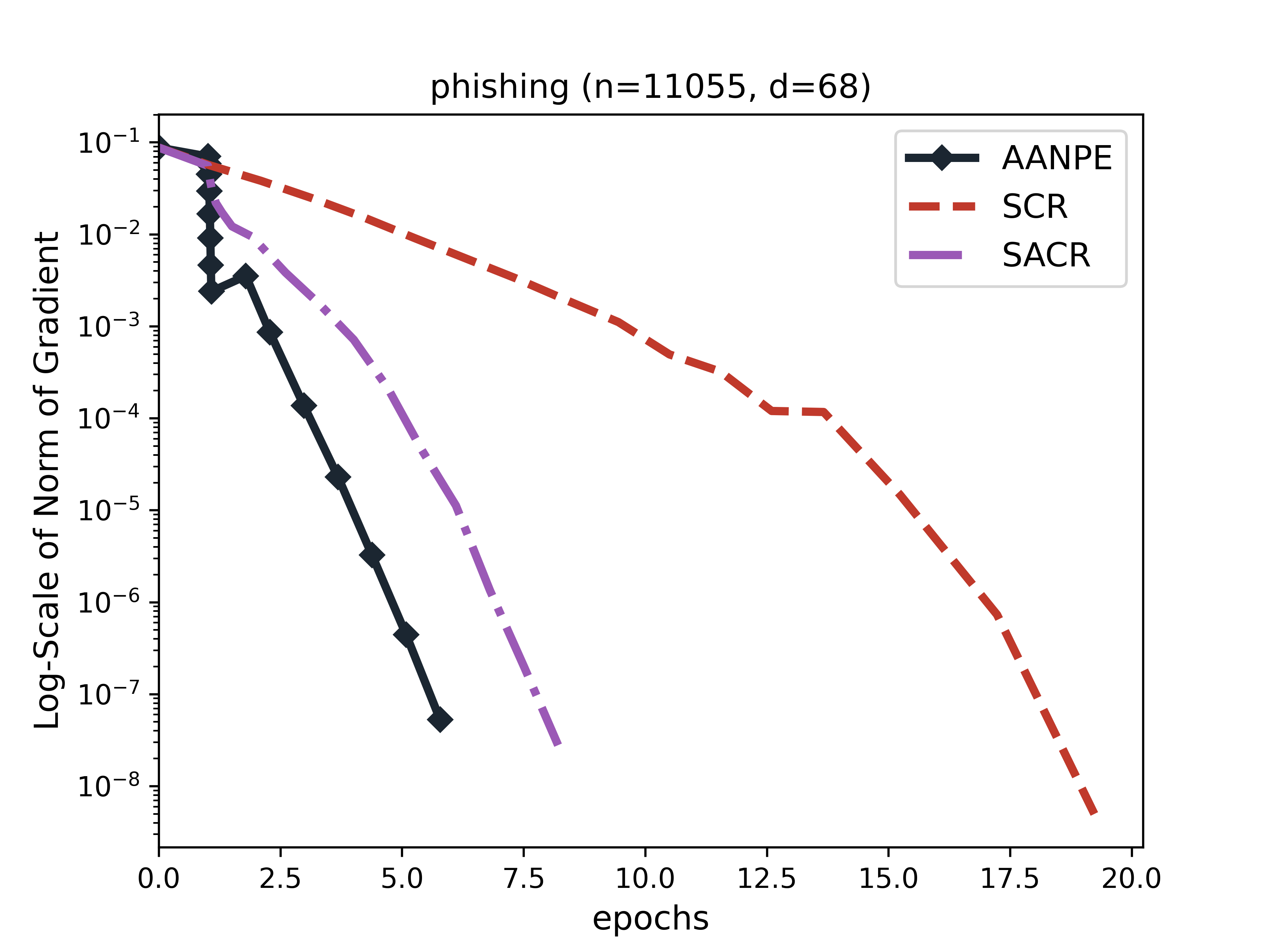}
  \end{minipage}
  \begin{minipage}[b]{0.49\textwidth}
    \includegraphics[width=\textwidth]{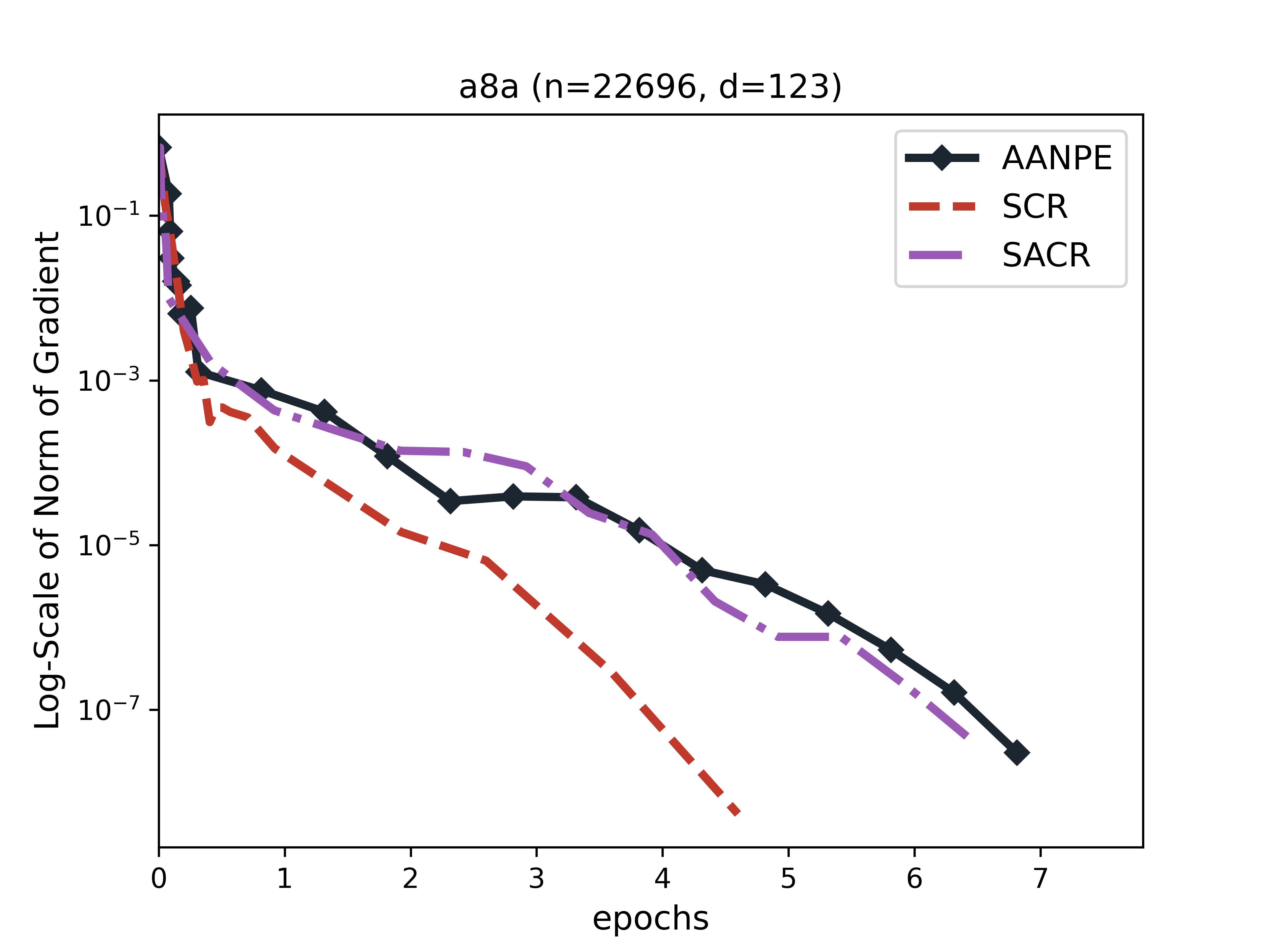}
  \end{minipage}
    \begin{minipage}[b]{0.49\textwidth}
    \includegraphics[width=\textwidth]{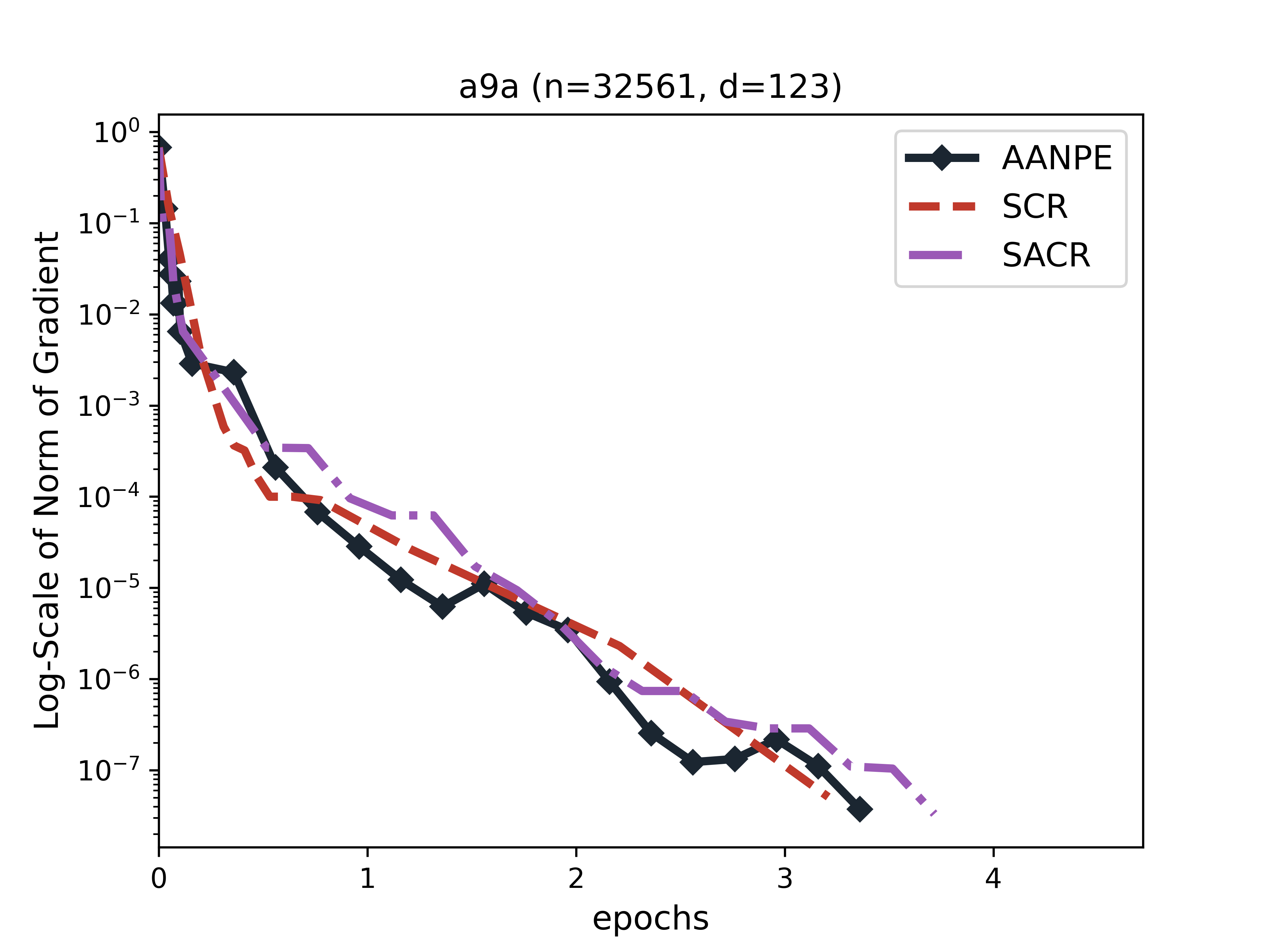}
  \end{minipage}
  \begin{minipage}[b]{0.49\textwidth}
    \includegraphics[width=\textwidth]{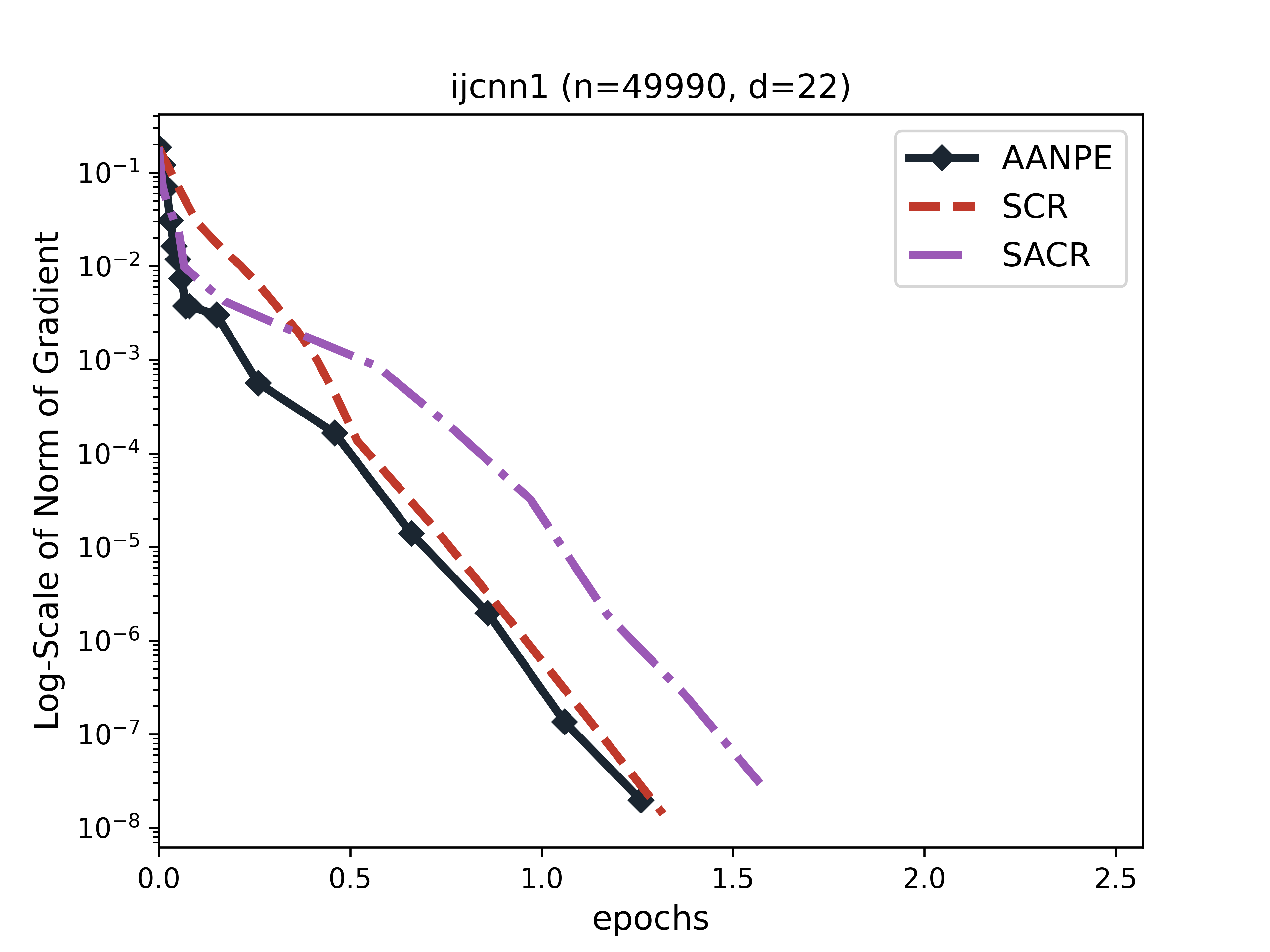}
  \end{minipage}
  \begin{minipage}[b]{0.49\textwidth}
    \includegraphics[width=\textwidth]{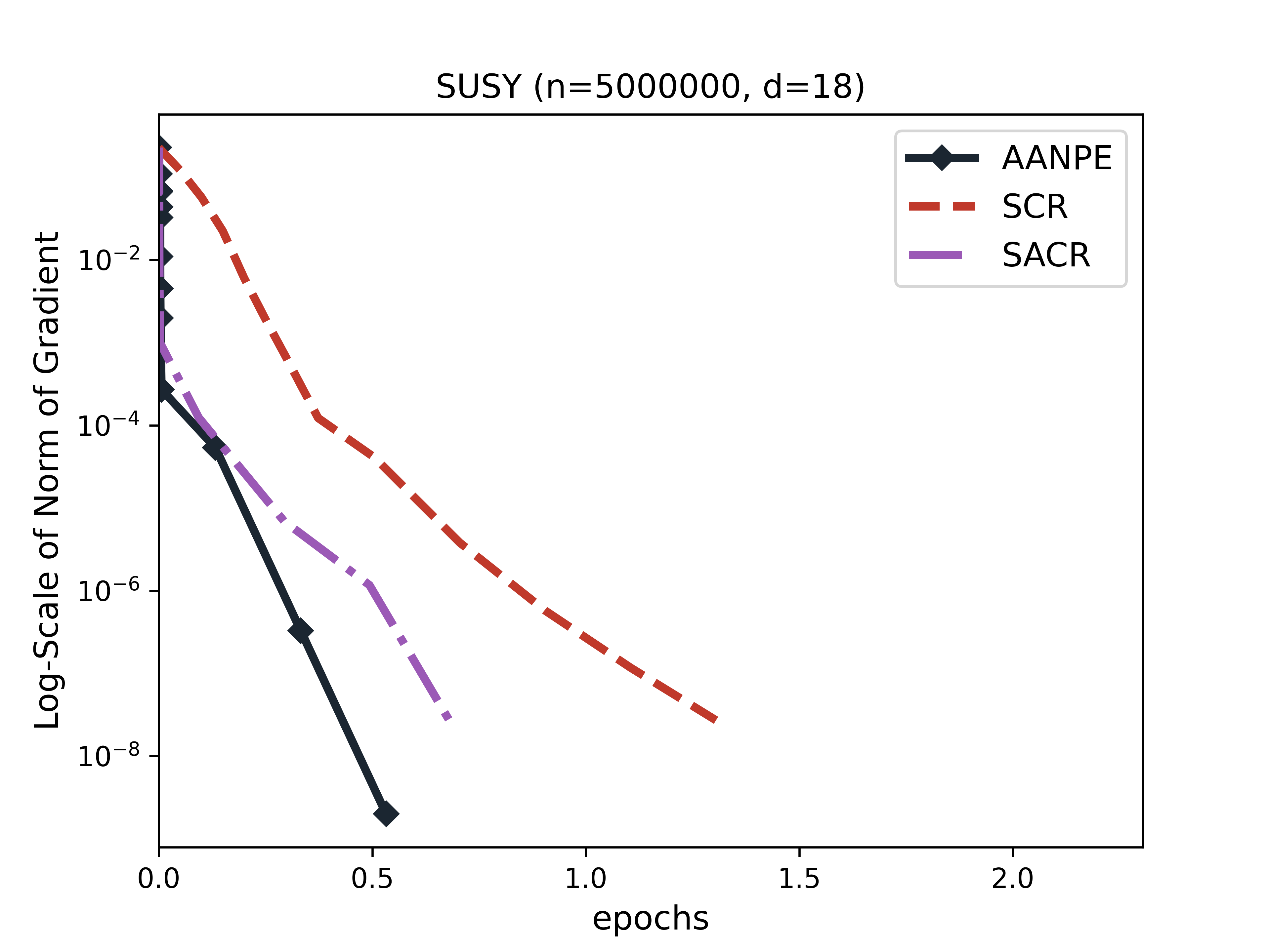}
  \end{minipage}
  \begin{minipage}[b]{0.49\textwidth}
    \includegraphics[width=\textwidth]{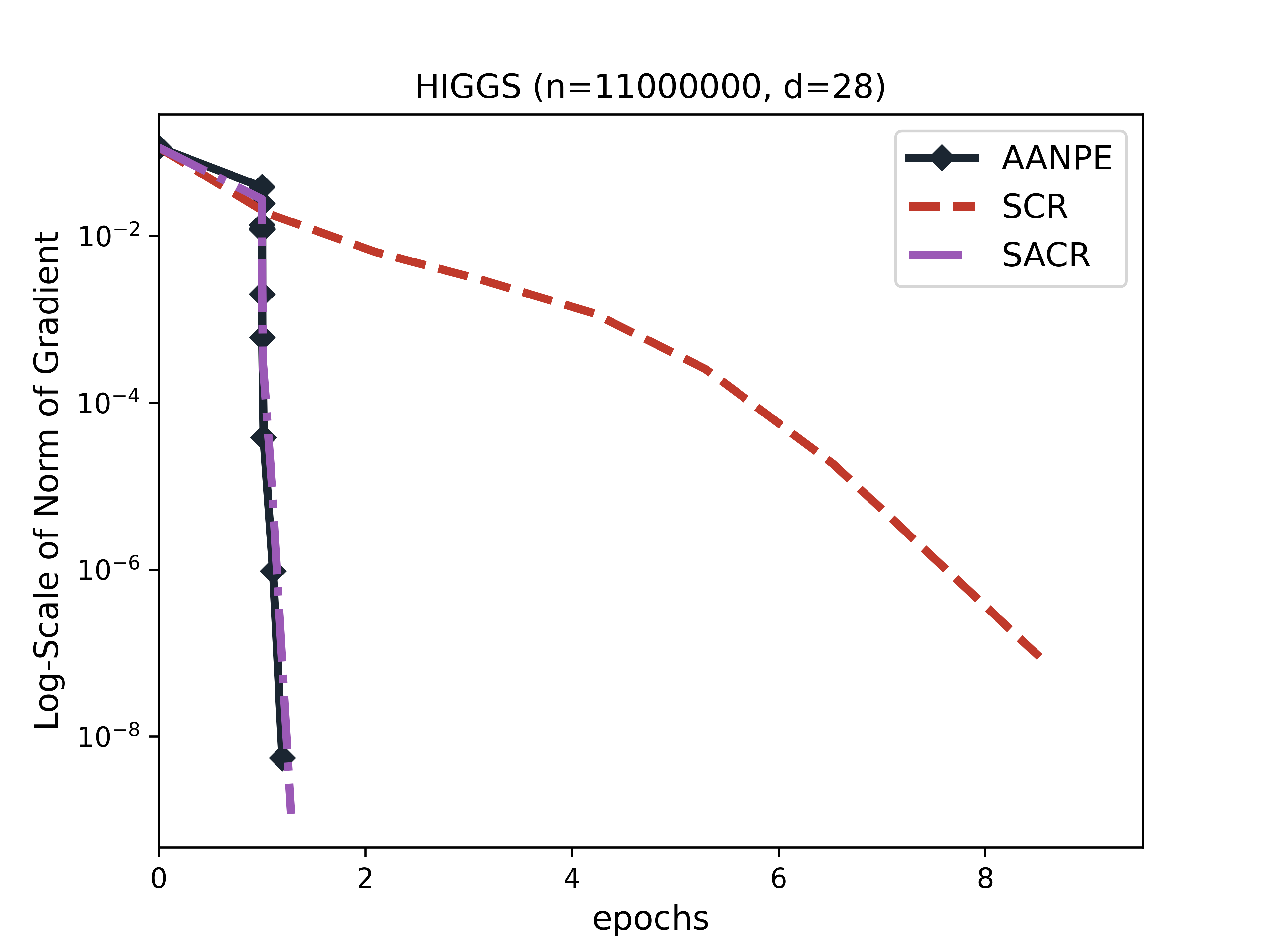}
  \end{minipage}
  \caption{
    Log-scale norm of gradient v.s. No. epochs.
  }
  \label{fig.random}
\end{figure}

We show the numerical results on the data sets in \autoref{fig.random}, where the log-scaled norm of gradient v.s. the number of epochs is shown. We can see from the picture that the three algorithms have similar convergence behavior and the performance of our algorithm is comparable to that of other two sub-sampled second-order methods.

\clearpage
\bibliographystyle{plainnat}


\appendix
\section{Proof of Technical Lemmas}
\subsection*{Proof to \autoref{prop:inexact approximate equation}}
{\it Proof}
We first note that \eqref{eq.bound v} and \eqref{eq.bound epsilon} can be derived by applying triangle inequality to \eqref{eq.inexact approxiamte equation lambda}. Thus it remains to prove \eqref{eq.inexact approxiamte equation v} and \eqref{eq.inexact approxiamte equation lambda}.

From the definition, we can easily verify that \eqref{eq.inexact approxiamte equation v} holds. As for \eqref{eq.inexact approxiamte equation lambda}, we can see that
\begin{equation}\label{eq.u-v}
    \begin{aligned}
        \|u-v\| & = \| \mathcal{G}(y)-\mathcal{G}_{x,\delta}(y)\|\\
        &=\|\mathcal{G}(y)-\mathcal{G}(x)-\mathcal{P}(x)(y-x)+\mathcal{G}'(x)(y-x)-\mathcal{G}'(x)(y-x)\|\\
        &\leq \|\mathcal{G}(y)-\mathcal{G}(x)-\mathcal{G}'(x)(y-x)\| + \|\left(\mathcal{P}\left(x\right)-\mathcal{G}'\left(x\right)\right)(y-x)\|\\
        &\leq \frac{L_2}{2}\|y-x\|^2 +\delta\|y-x\|.
    \end{aligned}
\end{equation}
Where the last line can be derived by the Lipschitz continuity of $\mathcal{G}'$, thus
\begin{equation*}
    \begin{aligned}
        \|\lambda v+y-x\|^2+2\lambda \epsilon &= \|\lambda u+y-x +\lambda v-\lambda u\|^2+2\lambda\epsilon\\
        &\leq \left ( \sqrt{\left \|\lambda u+y-x\right \|^2+2\lambda \epsilon}+\lambda \left \|u-x\right\| \right )^2\\
        &\leq \left ( \hat{\sigma}\left \|y-x\right \| +\frac{L_2}{2}\left \| y-x\right \| ^2 +\delta \left \|y-x\right \|\right )^2\\
        &=\left (\hat{\sigma}+\frac{\lambda L_2}{2}\left \|y-x\right \|+\lambda\delta\right )\|y-x\|^2.
    \end{aligned}
\end{equation*}
The second line can be verified by expanding the squared term on the right-hand side, the third line is from \autoref{def.approximate newton solution mmo}, and the last line is from \eqref{eq.u-v}.
\qed

\subsection*{Proof to \autoref{yuan7.5}}
{\it Proof}
For simplicity, let $y=y_\mathcal{B}(\lambda;x)$ and $\tilde{y}=y_{\tilde{\mathcal{B}}}(\lambda;\tilde{x})$, then there exist $v \in \mathcal{B}(y)$ and $\tilde{v} \in \tilde{\mathcal{B}}(\tilde{y})$ such that 
\begin{equation}\label{eq.two optcond}
    \lambda v+y-x=0, \ \lambda \tilde{v}+\tilde{y}-\tilde{x}=0.
\end{equation}
As a consequence,
$$\varphi_\mathcal{B}(\lambda;x)=\lambda^2\Vert v \Vert, \ \varphi_{\tilde{\mathcal{B}}}(\lambda;\tilde{x})=\lambda^2\Vert \tilde{v} \Vert.$$
Let $u:=v+\mathcal{G}_{\tilde{x},\delta}(y)-\mathcal{G}_{x,\delta}(y),$ note the that $v \in \mathcal{B}(y)$, then we know
\begin{equation*}
u\in \mathcal{B}(y)+\mathcal{G}_{\tilde{x},\delta}(y)-\mathcal{G}_{x,\delta}(y)=\mathcal{T}_{x,\delta}(y)+\mathcal{G}_{\tilde{x},\delta}(y)-\mathcal{G}_{x,\delta}(y)=\mathcal{G}_{\tilde{x},\delta}(y)+\mathcal{H}(y)=\tilde{\mathcal{B}}(y)
\end{equation*}
and
\begin{equation}\label{eq.error u}
\lambda u+y-\tilde{x}=\lambda v+y-x+(\tilde{x}-x)+\lambda(u-v)=(\tilde{x}-x)+\lambda(u-v).
\end{equation}
Combine \eqref{eq.error u} with \eqref{eq.two optcond} we conclude that 
\begin{equation}\label{eq.error u2}
\lambda(u-\tilde{v})+(y-\tilde{y})=(\tilde{x}-x)+\lambda(u-v).
\end{equation}
Since $u \in \tilde{\mathcal{B}}(y)$ and $\tilde{v} \in \tilde{\mathcal{B}}(\tilde{y})$, it follows from the monotonicity of $\tilde{\mathcal{B}}$ that $\langle u-\tilde{v},y-\tilde{y}\rangle \geq 0$, which together with \eqref{eq.error u2} and the triangle inequality for norms implies that
\begin{equation*}
\lambda \Vert u-\tilde{v} \Vert \le \Vert \tilde{x}-x \Vert +\lambda \Vert u-v \Vert
\end{equation*}
and hence that
\begin{equation*}
\lambda \Vert v-\tilde{v} \Vert \le \Vert \tilde{x}-x \Vert +2\lambda \Vert u-v \Vert.
\end{equation*}
This implies
\begin{equation*}
\vert \varphi_\mathcal{B}(\lambda;x)-\varphi_{\tilde{\mathcal{B}}}(\lambda;\Tilde{x}) \vert=\lambda^2 \left \vert \left \Vert v \right\Vert-\left \Vert \tilde{v}\right \Vert \right \vert \le \lambda^2 \left \Vert v-\tilde{v}\right \Vert \le \lambda \left ( \left \Vert \tilde{x}-x \right \Vert +2\lambda \left \Vert u-v \right \Vert \right).
\end{equation*}
Now, use the definition of $u$, we have
\begin{align*}
u-v&=\mathcal{G}_{\tilde{x},\delta}(y)-\mathcal{G}_{x,\delta}(y)=\mathcal{G}(\tilde{x})+\mathcal{P}(\tilde{x})(y-\tilde{x})-\left (\mathcal{G}(x)+\mathcal{P}(x)(y-x)\right )\\
&=\mathcal{G}(\tilde{x})+\mathcal{P}(\tilde{x})(x-\tilde{x})-\mathcal{G}(x)+\left (\mathcal{P}(\tilde{x})-\mathcal{P}(x)\right )(y-x)
\end{align*}
and hence
\begin{align*}
\lambda \Vert u-v \Vert & \le \lambda \Vert \mathcal{G}(\tilde{x})+\mathcal{P}(\tilde{x})(x-\tilde{x})-\mathcal{G}(x) \Vert+\lambda \Vert \mathcal{P}(\tilde{x})-\mathcal{P}(x) \Vert \Vert y-x \Vert\\
& \le \lambda \Vert \mathcal{G}(\tilde{x})+\mathcal{G}'(\tilde{x})(x-\tilde{x})-\mathcal{G}(x) \Vert + \lambda \Vert \mathcal{P}(\tilde{x})-\mathcal{G}'(\tilde{x})\vert \Vert x-\tilde{x}   \Vert +\lambda \Vert \mathcal{P}(\tilde{x})-\mathcal{P}(x) \Vert \Vert y-x \Vert \\
& \le \frac{\lambda L_2}{2} \Vert \tilde{x}-x \Vert^2+\lambda \delta \Vert \tilde{x}-x \Vert+\lambda L_2\Vert \tilde{x}-x \Vert  \Vert y-x \Vert\\
& \le \frac{\lambda L_2}{2} \Vert \tilde{x}-x \Vert^2+\lambda \delta \Vert \tilde{x}-x \Vert+\lambda L_2\Vert \tilde{x}-x \Vert  \Vert y-x \Vert \\
&=\frac{\lambda L_2}{2} \Vert \tilde{x}-x \Vert^2+\lambda \delta \Vert \tilde{x}-x \Vert+ L_2\Vert \tilde{x}-x \Vert  \varphi_{\mathcal{B}}(\lambda;x)\\
& \le  \frac{\lambda L_2}{2} \Vert \tilde{x}-x \Vert^2+C\Vert \tilde{x}-x \Vert+ L_2\Vert \tilde{x}-x \Vert  \varphi_{\mathcal{B}}(\lambda;x)
\end{align*}
Now we can conclude that
\begin{equation*}
\vert \varphi_\mathcal{B}(\lambda;x)-\varphi_{\tilde{\mathcal{B}}}(\lambda;\tilde{x})\vert \le (1+2C)\lambda \Vert \tilde{x}-x \Vert+L_2\lambda^2 \Vert \tilde{x}-x \Vert^2+2L_2\lambda \Vert \tilde{x}-x \Vert \varphi_\mathcal{B}(\lambda;x).
\end{equation*}
This inequality and the symmetric one obtained by interchanging $x$ and $\tilde{x}$ in the latter relation then imply \eqref{varphi-varphi}.
\qed

\subsection*{Proof to \autoref{prop.right bracketing point}}
{\it Proof}
Suppose $\lambda \geq \max \left\{ \sqrt{\frac{\alpha}{\bar{v}}\left(1+\hat{\sigma}+C+\frac{L_2\alpha}{2}\right)},  \left(\frac{\hat{\sigma}^2 \alpha^2}{2\bar{\epsilon}}\right)^{\frac{1}{3}} \right\}$ and $\lambda \Vert y-x \Vert \le \alpha$. From proposition~\ref{prop:inexact approximate equation} we know that 
\begin{equation*}
\Vert v \Vert \le \frac{1}{\lambda}\left(1+\hat{\sigma}+C+\frac{L_2\lambda}{2} \Vert y-x \Vert\right)\Vert y-x \Vert \le \left(1+\hat{\sigma}+C+\frac{L_2\alpha}{2}\right) \frac{\alpha}{\lambda^2} \le \bar{v}
\end{equation*}
and
\begin{equation*}
\epsilon \le \frac{\hat{\sigma}^2}{2\lambda} \Vert y-x \Vert^2 \le \frac{\hat{\sigma}^2 \alpha^2}{2\lambda^3} \le \bar{\epsilon}.
\end{equation*}
This contradicts \eqref{eq.relation v rho}.
\qed

\subsection*{Proof to \autoref{prop:min}}
{\it Proof}
First, we note that \eqref{answer} follows immediately from \eqref{eq.lmpz} and \eqref{eq.lmmz}. Let $\mathcal{B}_-:=\mathcal{T}_{x_-^0,\delta}$ and $\mathcal{B}_+:=\mathcal{T}_{x_+^0,\delta}$. Since a $(\hat{\sigma},\delta)$-approximate Newton solution at $(\lambda_-^0,x_-^0)$ is a $\hat{\sigma}$-approximate solution with $\mathcal{B}=\mathcal{B}_-$, it follows from Proposition~\ref{yinyong2} with $\mathcal{B}=\mathcal{B}_-$ that
\begin{equation*}
\lambda_-^0 \Vert y_-^0-x_-^0 \Vert \le \frac{\varphi_{\mathcal{B}_-}(\lambda_-^0;x_-^0)}{1-\hat{\sigma}} \le \frac{\lambda_-^0 \varphi_{\mathcal{B}_-}(\lambda_+^0;x_-^0)}{(1-\hat{\sigma})\lambda_+^0},
\end{equation*}
where the last inequality is due to \eqref{answer} and Proposition~\ref{yinyong1}(b) with $\mathcal{B}=\mathcal{B}_-, \tilde{\lambda}=\lambda_-^0$, and $\lambda=\lambda_+^0$. Also Proposition~\ref{yuan7.5} with $\lambda=\lambda_+^0,x=x_-^0$, and $\tilde{x}=x_+^0$ and the definition of $\theta_+^0$ implies that
\begin{align*}
\varphi_{\mathcal{B}_-}(\lambda_+^0;x_-^0) &\le (2L_2\theta_+^0+1)\varphi_{\mathcal{B}_+}(\lambda_+^0;x_+^0)+(1+2C)\theta_+^0+L_2(\theta_+^0)^2\\
&\le (2L_2\theta_+^0+1)(1+\hat{\sigma})\lambda_+^0 \Vert y_+^0-x_+^0 \Vert+(1+2C)\theta_+^0+L_2(\theta_+^0)^2\\
&\le \frac{\alpha (1-\hat{\sigma})\lambda_+^0}{\lambda_-^0},
\end{align*}
where the last inequality follows from the definition of $\lambda_-^0$. Combining the above two inequalities, we then conclude that \eqref{conclusion} holds.
\qed

\subsection*{Proof to \autoref{coro.left bracketing point}}
{\it Proof}
From \eqref{eq.settolerance} and \autoref{prop.right bracketing point} we can easily know that the first claim holds.
From \eqref{eq.left bracketing point} we know
\begin{equation*}
\lambda_{-}^0  \leq \frac{\alpha_- (1-\hat{\sigma})\lambda_+^0}{(1+\hat{\sigma})(1+2L_2\theta_+^0)\lambda_+^0 \Vert y_+^0-x_+^0 \Vert+(1+2C)\theta_+^0+L_2(\theta_+^0)^2},
\end{equation*}
so in the view of Proposition~\ref{prop:min} we have
\begin{equation*}
\lambda_-^0 \Vert y_-^0-x_-^0\Vert \leq \alpha_-.
\end{equation*}
\qed

\subsection*{Proof to \autoref{lem.length while loop}}
{\it Proof}
    Note that as in \autoref{prop.right bracketing point}, when $\lambda_{k+1}\geq \Lambda\left(\frac{2\sigma_l}{L_2}\right)$, the while loop will terminate. When the length of the while loop approaches $\log_{\gamma}\frac{\delta_{\max}\Lambda\left(\frac{2\sigma_l}{L_2}\right)}{C}$, we have $\lambda_{k+1}=\frac{C}{\delta_k}\geq \Lambda\left(\frac{2\sigma_l}{L_2}\right)$.
\qed

\subsection*{Proof to \autoref{prop:change}}
{\it Proof}
Let $\mathcal{B}$ denote the exact first-order approximate of $\mathcal{T}$ at $x$, then from Monteiro and Svaiter \cite[Theorem 7.11]{monteiro2013accelerated} we know that:
\begin{equation}\label{oldthm}
\varphi_\mathcal{B}(\lambda;x)\le \lambda \Vert x-x_*\Vert+\lambda^2 L_2 \Vert x-x_* \Vert^2.
\end{equation}
Since
\begin{equation*}
\varphi_\mathcal{B}(\lambda;x)=\lambda \Vert x-(I+\lambda \mathcal{B})^{-1}(x) \Vert,
\end{equation*}
\begin{equation*}
\varphi_{\mathcal{T}_{x,\delta}}(\lambda;x)=\lambda \Vert x-(I+\lambda \mathcal{T}_{x,\delta})^{-1}(x) \Vert.
\end{equation*}
So we let $y_1=(I+\lambda \mathcal{B})^{-1}(x)$ and $y_2=(I+\lambda \mathcal{T}_{x,\delta})^{-1}(x)$, then we have
\begin{align}
\label{eq.boundinexact_exact}
\varphi_\mathcal{B}(\lambda;x)&=\lambda \Vert x-y_1 \Vert,\\
\label{eq.boundinexact_inexact}
\varphi_{\mathcal{T}_{x,\delta}}(\lambda;x)&=\lambda \Vert x-y_2 \Vert,
\end{align}
and
\begin{align}
\label{eq.boundinexact_bridge1}
y_1+\lambda u=x, \ u = \mathcal{G}_x(y)+u_h, \ u_h\in \mathcal{H}(y_1),\\
\label{eq.boundinexact_bridge2}
y_2+\lambda v=x, \ v = \mathcal{G}_{x,\delta}(y)+v_h, \ v_h\in \mathcal{H}(y_2).
\end{align}
$\mathcal{G}_x(y):= \mathcal{G}(x)+\mathcal{G}'(x)(y-x)$. Subtracting the equations \eqref{eq.boundinexact_bridge1} and \eqref{eq.boundinexact_bridge2},
\begin{align*}
0 &=y_1-y_2+\lambda (u-v)\\
&=y_1-y_2+\lambda \mathcal{G}'(x)(y_1-x)+\lambda u_h-\lambda \mathcal{P}(x)(y_2-x)-\lambda v_h\\
&=(y_1-y_2)+\lambda \mathcal{G}'(x)y_1-\lambda \mathcal{P}(x)y_2+\lambda (\mathcal{P}(x)-\mathcal{G}'(x))x+\lambda(u_h-v_h)\\
&=(y_1-y_2)+\lambda \mathcal{G}'(x)(y_1-y_2)-\lambda (\mathcal{P}(x)-\mathcal{G}'(x))y_2+\lambda (\mathcal{P}(x)-\mathcal{G}'(x))x+\lambda(u_h-v_h).
\end{align*}
Then we know
\begin{equation*}
(1+\lambda \mathcal{G}'(x))(y_1-y_2)=(\mathcal{P}(x)-\mathcal{G}'(x))\lambda (y_2-x)-\lambda(u_h-v_h).
\end{equation*}
Note that $\mathcal{G}$ is monotone and differentiable, as a consequence, 
\begin{equation}
\label{eq.bound_y1y2}
\begin{aligned}
\Vert y_1-y_2 \Vert &\le \lambda\Vert \mathcal{P}(x)-\mathcal{G}'(x) \Vert \Vert  y_2-x \Vert + \lambda\Vert u_h-v_h \Vert\\
& \le \delta \lambda\Vert  y_2-x \Vert +2L'\lambda\\
&=\delta \varphi_{\mathcal{T}_x}(\lambda;x)+2L'\lambda,
\end{aligned}    
\end{equation}
the first inequality is from that $\mathcal{G}$ is maximal monotone, the second inequality comes from the boundedness of $\mathcal{H}$.
From the above inequation, \eqref{eq.boundinexact_inexact} and \eqref{oldthm} we know
\begin{equation}
\label{eq.bound_inexact}
\begin{aligned}
\Vert x-(I+\lambda \mathcal{T}_x)^{-1}(x) \Vert & = \| x-y_2\| \\
& \le \| x-y_1 \|+\| y_1-y_2\| \\
&\le \delta \varphi_{\mathcal{T}_x}(\lambda;x)+2L'\lambda +\Vert x-x_*\Vert+\lambda L_2 \Vert x-x_* \Vert^2,
\end{aligned}
\end{equation}
note that $\lambda\delta \le C$, \eqref{eq.bound_inexact} implies
\begin{equation}\label{eq.bound_inexactfinal}
\varphi_{\mathcal{T}_x}(\lambda;x) \le C\varphi_{\mathcal{T}_x}(\lambda;x)+2L'\lambda^2+\lambda \Vert x-x_*\Vert+\lambda^2 L_2 \Vert x-x_* \Vert^2.
\end{equation}
Then the proof is finished.
\qed

\subsection*{Proof to \autoref{lemma:h}}
{\it Proof}
For any $v\in \partial_\epsilon h(x)$, let $y_\alpha = x+ \alpha v$, from \eqref{epsilon-subdifferential}, we have
\begin{equation*}
    \begin{aligned}
        \alpha \|v\|^2 &= \langle y-x,v\rangle\\
        &\leq h(y)-h(x)+\epsilon\\
        &\leq L'\|y-x\|+\epsilon\\
        &=\alpha L'\|v\|+\epsilon.
    \end{aligned}
\end{equation*}
The second last line comes from the Lipschitz continuity of 
$h$. Let $\alpha \to +\infty$ we have $\|v\| \le L'$.
\qed

\subsection*{Proof to \autoref{cor:uniform}}
{\it Proof}
Since sample size
\begin{equation*}
\left|\SCal\right| \geq \frac{16L_1^2}{\kappa^2} \log\left(\frac{2Nd}{\delta'}\right),
\end{equation*}
as in \cite{xu2016sub}, for each $i=1,2,\dotsb,N$ we have 
\begin{equation*}
\Prob\left(\left\| \Tilde{H}(x_i) - \nabla^2 g(x_i)\right\| \geq \kappa \right) < \frac{\delta'}{N}.
\end{equation*}
Then the corollary is from the subadditivity of the probability.
\qed
\end{document}